

\overfullrule=0pt
\def\title{Vector bundles on elliptic curve and Sklyanin algebras}

\def\author{B. L. Feigin and A. V. Odesskii}

\magnification=\magstep1
\baselineskip=18pt plus 3pt
\font\twelvebf=cmbx12

\font\bigmi=cmmi10 scaled \magstep1
\font\bigsy=cmsy10 scaled \magstep1
\font\sbigmi=cmmi10
\font\sbigsy=cmsy10
\font\ssbigmi=cmmi8
\font\ssbigsy=cmsy8
\def\bigbffont{\textfont0=\twelvebf
     \textfont1=\bigmi \scriptfont1=\sbigmi \scriptscriptfont1=\ssbigmi
     \textfont2=\bigsy \scriptfont2=\sbigsy \scriptscriptfont2=\ssbigsy
     \twelvebf}
\vglue1.5cm

\centerline{\bigbffont
\vbox{
\halign{\hfil # \hfil\cr
\title\crcr}}}

\bigskip\bigskip

\centerline{\tenrm
\vbox{
\halign{\hfil#\hfil\cr
\author\crcr}}}

\bigskip\bigskip

\def\eq#1\endeq
{$$\eqalignno{#1}$$}

\def\leq#1\endeq
{$$\leqalignno{#1}$$}

\def\proof#1\par{\medbreak\noindent{\it Proof.\/}\quad{#1}\par
    \ifdim\lastskip<\medskipamount\removelastskip\penalty55\medskip\fi}

\def\remark#1\par{\medbreak\noindent{\it Remark.}\quad#1\par
    \ifdim\lastskip<\medskipamount\removelastskip\penalty55\medskip\fi}

\outer\def\Beginsection#1#2\par{\vskip0pt plus.3\vsize\penalty-250
  \vskip0pt plus-.3\vsize\bigskip\vskip\parskip
  \message{#1}\noindent\hangindent=#2 \hangafter=1
  {\bf#1}\smallskip\noindent}

\catcode`\@=11

\font\tenmsa=msam10
\font\sevenmsa=msam7
\font\fivemsa=msam5
\font\tenmsb=msbm10
\font\sevenmsb=msbm7
\font\fivemsb=msbm5
\newfam\msafam
\newfam\msbfam
\textfont\msafam=\tenmsa  \scriptfont\msafam=\sevenmsa
  \scriptscriptfont\msafam=\fivemsa
\textfont\msbfam=\tenmsb  \scriptfont\msbfam=\sevenmsb
  \scriptscriptfont\msbfam=\fivemsb

\def\hexnumber@#1{\ifnum#1<10 \number#1\else
 \ifnum#1=10 A\else\ifnum#1=11 B\else\ifnum#1=12 C\else
 \ifnum#1=13 D\else\ifnum#1=14 E\else\ifnum#1=15 F\fi\fi\fi\fi\fi\fi\fi}

\def\msa@{\hexnumber@\msafam}
\def\msb@{\hexnumber@\msbfam}
\mathchardef\boxdot="2\msa@00
\mathchardef\boxplus="2\msa@01
\mathchardef\boxtimes="2\msa@02
\mathchardef\square="0\msa@03
\mathchardef\blacksquare="0\msa@04
\mathchardef\centerdot="2\msa@05
\mathchardef\lozenge="0\msa@06
\mathchardef\blacklozenge="0\msa@07
\mathchardef\circlearrowright="3\msa@08
\mathchardef\circlearrowleft="3\msa@09
\mathchardef\rightleftharpoons="3\msa@0A
\mathchardef\leftrightharpoons="3\msa@0B
\mathchardef\boxminus="2\msa@0C
\mathchardef\Vdash="3\msa@0D
\mathchardef\Vvdash="3\msa@0E
\mathchardef\vDash="3\msa@0F
\mathchardef\twoheadrightarrow="3\msa@10
\mathchardef\twoheadleftarrow="3\msa@11
\mathchardef\leftleftarrows="3\msa@12
\mathchardef\rightrightarrows="3\msa@13
\mathchardef\upuparrows="3\msa@14
\mathchardef\downdownarrows="3\msa@15
\mathchardef\upharpoonright="3\msa@16

\mathchardef\downharpoonright="3\msa@17
\mathchardef\upharpoonleft="3\msa@18
\mathchardef\downharpoonleft="3\msa@19
\mathchardef\rightarrowtail="3\msa@1A
\mathchardef\leftarrowtail="3\msa@1B
\mathchardef\leftrightarrows="3\msa@1C
\mathchardef\rightleftarrows="3\msa@1D
\mathchardef\Lsh="3\msa@1E
\mathchardef\Rsh="3\msa@1F
\mathchardef\rightsquigarrow="3\msa@20
\mathchardef\leftrightsquigarrow="3\msa@21
\mathchardef\looparrowleft="3\msa@22
\mathchardef\looparrowright="3\msa@23
\mathchardef\circeq="3\msa@24
\mathchardef\succsim="3\msa@25
\mathchardef\gtrsim="3\msa@26
\mathchardef\gtrapprox="3\msa@27
\mathchardef\multimap="3\msa@28
\mathchardef\therefore="3\msa@29
\mathchardef\because="3\msa@2A
\mathchardef\doteqdot="3\msa@2B

\mathchardef\triangleq="3\msa@2C
\mathchardef\precsim="3\msa@2D
\mathchardef\lesssim="3\msa@2E
\mathchardef\lessapprox="3\msa@2F
\mathchardef\eqslantless="3\msa@30
\mathchardef\eqslantgtr="3\msa@31
\mathchardef\curlyeqprec="3\msa@32
\mathchardef\curlyeqsucc="3\msa@33
\mathchardef\preccurlyeq="3\msa@34
\mathchardef\leqq="3\msa@35
\mathchardef\leqslant="3\msa@36
\mathchardef\lessgtr="3\msa@37
\mathchardef\backprime="0\msa@38
\mathchardef\risingdotseq="3\msa@3A
\mathchardef\fallingdotseq="3\msa@3B
\mathchardef\succcurlyeq="3\msa@3C
\mathchardef\geqq="3\msa@3D
\mathchardef\geqslant="3\msa@3E
\mathchardef\gtrless="3\msa@3F
\mathchardef\sqsubset="3\msa@40
\mathchardef\sqsupset="3\msa@41
\mathchardef\trianglerighteq="3\msa@44
\mathchardef\trianglelefteq="3\msa@45
\mathchardef\bigstar="0\msa@46
\mathchardef\between="3\msa@47
\mathchardef\blacktriangledown="0\msa@48
\mathchardef\blacktriangleright="3\msa@49
\mathchardef\blacktriangleleft="3\msa@4A
\mathchardef\blacktriangle="0\msa@4E
\mathchardef\triangledown="0\msa@4F
\mathchardef\eqcirc="3\msa@50
\mathchardef\lesseqgtr="3\msa@51
\mathchardef\gtreqless="3\msa@52
\mathchardef\lesseqqgtr="3\msa@53
\mathchardef\gtreqqless="3\msa@54
\mathchardef\Rrightarrow="3\msa@56
\mathchardef\Lleftarrow="3\msa@57
\mathchardef\veebar="2\msa@59
\mathchardef\barwedge="2\msa@5A
\mathchardef\doublebarwedge="2\msa@5B
\mathchardef\angle="0\msa@5C
\mathchardef\measuredangle="0\msa@5D
\mathchardef\sphericalangle="0\msa@5E
\mathchardef\varpropto="3\msa@5F
\mathchardef\smallsmile="3\msa@60
\mathchardef\smallfrown="3\msa@61
\mathchardef\Subset="3\msa@62
\mathchardef\Supset="3\msa@63
\mathchardef\Cup="2\msa@64

\mathchardef\Cap="2\msa@65

\mathchardef\curlywedge="2\msa@66
\mathchardef\curlyvee="2\msa@67
\mathchardef\leftthreetimes="2\msa@68
\mathchardef\rightthreetimes="2\msa@69
\mathchardef\subseteqq="3\msa@6A
\mathchardef\supseteqq="3\msa@6B
\mathchardef\bumpeq="3\msa@6C
\mathchardef\Bumpeq="3\msa@6D
\mathchardef\lll="3\msa@6E

\mathchardef\ggg="3\msa@6F

\mathchardef\circledS="0\msa@73
\mathchardef\pitchfork="3\msa@74
\mathchardef\dotplus="2\msa@75
\mathchardef\backsim="3\msa@76
\mathchardef\backsimeq="3\msa@77
\mathchardef\complement="0\msa@7B
\mathchardef\intercal="2\msa@7C
\mathchardef\circledcirc="2\msa@7D
\mathchardef\circledast="2\msa@7E
\mathchardef\circleddash="2\msa@7F
\def\ulcorner{\delimiter"4\msa@70\msa@70 }
\def\urcorner{\delimiter"5\msa@71\msa@71 }
\def\llcorner{\delimiter"4\msa@78\msa@78 }
\def\lrcorner{\delimiter"5\msa@79\msa@79 }
\def\yen{\mathhexbox\msa@55 }
\def\checkmark{\mathhexbox\msa@58 }
\def\circledR{\mathhexbox\msa@72 }
\def\maltese{\mathhexbox\msa@7A }
\mathchardef\lvertneqq="3\msb@00
\mathchardef\gvertneqq="3\msb@01
\mathchardef\nleq="3\msb@02
\mathchardef\ngeq="3\msb@03
\mathchardef\nless="3\msb@04
\mathchardef\ngtr="3\msb@05
\mathchardef\nprec="3\msb@06
\mathchardef\nsucc="3\msb@07
\mathchardef\lneqq="3\msb@08
\mathchardef\gneqq="3\msb@09
\mathchardef\nleqslant="3\msb@0A
\mathchardef\ngeqslant="3\msb@0B
\mathchardef\lneq="3\msb@0C
\mathchardef\gneq="3\msb@0D
\mathchardef\npreceq="3\msb@0E
\mathchardef\nsucceq="3\msb@0F
\mathchardef\precnsim="3\msb@10
\mathchardef\succnsim="3\msb@11
\mathchardef\lnsim="3\msb@12
\mathchardef\gnsim="3\msb@13
\mathchardef\nleqq="3\msb@14
\mathchardef\ngeqq="3\msb@15
\mathchardef\precneqq="3\msb@16
\mathchardef\succneqq="3\msb@17
\mathchardef\precnapprox="3\msb@18
\mathchardef\succnapprox="3\msb@19
\mathchardef\lnapprox="3\msb@1A
\mathchardef\gnapprox="3\msb@1B
\mathchardef\nsim="3\msb@1C
\mathchardef\napprox="3\msb@1D
\mathchardef\nsubseteqq="3\msb@22
\mathchardef\nsupseteqq="3\msb@23
\mathchardef\subsetneqq="3\msb@24
\mathchardef\supsetneqq="3\msb@25
\mathchardef\subsetneq="3\msb@28
\mathchardef\supsetneq="3\msb@29
\mathchardef\nsubseteq="3\msb@2A
\mathchardef\nsupseteq="3\msb@2B
\mathchardef\nparallel="3\msb@2C
\mathchardef\nmid="3\msb@2D
\mathchardef\nshortmid="3\msb@2E
\mathchardef\nshortparallel="3\msb@2F
\mathchardef\nvdash="3\msb@30
\mathchardef\nVdash="3\msb@31
\mathchardef\nvDash="3\msb@32
\mathchardef\nVDash="3\msb@33
\mathchardef\ntrianglerighteq="3\msb@34
\mathchardef\ntrianglelefteq="3\msb@35
\mathchardef\ntriangleleft="3\msb@36
\mathchardef\ntriangleright="3\msb@37
\mathchardef\nleftarrow="3\msb@38
\mathchardef\nrightarrow="3\msb@39
\mathchardef\nLeftarrow="3\msb@3A
\mathchardef\nRightarrow="3\msb@3B
\mathchardef\nLeftrightarrow="3\msb@3C
\mathchardef\nleftrightarrow="3\msb@3D
\mathchardef\divideontimes="2\msb@3E
\mathchardef\varnothing="0\msb@3F
\mathchardef\nexists="0\msb@40
\mathchardef\mho="0\msb@66
\mathchardef\thorn="0\msb@67
\mathchardef\beth="0\msb@69
\mathchardef\gimel="0\msb@6A
\mathchardef\daleth="0\msb@6B
\mathchardef\lessdot="3\msb@6C
\mathchardef\gtrdot="3\msb@6D
\mathchardef\ltimes="2\msb@6E
\mathchardef\rtimes="2\msb@6F
\mathchardef\shortmid="3\msb@70
\mathchardef\shortparallel="3\msb@71
\mathchardef\smallsetminus="2\msb@72
\mathchardef\thicksim="3\msb@73
\mathchardef\thickapprox="3\msb@74
\mathchardef\approxeq="3\msb@75
\mathchardef\succapprox="3\msb@76
\mathchardef\precapprox="3\msb@77
\mathchardef\curvearrowleft="3\msb@78
\mathchardef\curvearrowright="3\msb@79
\mathchardef\digamma="0\msb@7A
\mathchardef\varkappa="0\msb@7B
\mathchardef\hslash="0\msb@7D
\mathchardef\hbar="0\msb@7E
\mathchardef\backepsilon="3\msb@7F
\def\Bbb{\ifmmode\let\next\Bbb@\else
 \def\next{\errmessage{Use \string\Bbb\space only in math mode}}\fi\next}
\def\Bbb@#1{{\Bbb@@{#1}}}
\def\Bbb@@#1{\fam\msbfam#1}

\catcode`\@=\active

\font\germ=eufm10
\def\goth#1{\hbox{\germ #1}}
\def\Ga{{\goth a}}

\def\GG{{\goth G}}

\def\BC{{\Bbb C}}
\def\BP{{\Bbb P}}
\def\BZ{{\Bbb Z}}

\def\CO{{\cal O}}
\def\CE{{\cal E}}
\def\cc{{\cal C}}

\def\Ch{{\rm Ch}}
\def\Ext{{\rm Ext}}
\def\Hom{{\rm Hom}}
\def\Mod{{\rm Mod}}

\def\al{\alpha}
\def\hrw{\hookrightarrow}
\def\rtw{\rightarrow}
\def\da{\Big\downarrow}

\def\maru#1{{\rm\ooalign{\hfil\lower.168ex\hbox{#1}%
\hfil\crcr\mathhexbox20D}}}
\def\ur#1{{\buildrel{\rm #1}\over\longrightarrow}}

\def\bk{{\bar k}}
\def\bn{{\bar n}}
\def\bs{{\bar S}}
\def\bM{{\bar M}}
\def\bal{{\bar \alpha}}
\def\bba{{\bar \beta}}
\def\bta{{\bar \tau}}

\def\apart{\noalign{\vskip8pt}}

\def\longrightrightarrow{\rlap{\lower1mm\hbox{$\longmapsto$}}
  \raise.7mm\hbox{$\longmapsto$}}
\def\yaya{\mathop{\longrightrightarrow}\limits}

\def\longrightrightmapsto{\rlap{\lower2mm\hbox{$\searrow$}}
  \raise2mm\hbox{$\nearrow$}}
\def\ayaya{\mathop{\longrightrightmapsto}\limits}

\beginsection
Introduction

In [4] we introduce the associative algebras
$Q_{n,k}(\CE,\tau)$.
Recall the definition.
These algebras are labeled by discrete parameters
$n,k$; $n,k$ are integers $n>k>0$
and $n$ and $k$ have not common divisors.
Then, $\CE$ is an elliptic curve and
$\tau$ is a point in $\CE$.
We identify $\CE$ with $\BC/\Gamma$,
where $\Gamma$ is a lattice.
\par
Algebra
$Q_{n,k}(\CE,\tau)$
is generated by $n$ generators
$\{x_i\}$,
$i\in\BZ/n\BZ$,
which satisfy the relations:
\eq
&\sum_{r\in\BZ/n\BZ}
{\theta_{j-i +
(k-1)r}(0)\over
\theta_{j-i-r}(-\tau)\theta_{kr}(\tau)}
x_{k(j-r)}
x_{k(i+r)}=0.\cr
\endeq
Here all indices belong to
$\BZ/n\BZ$, $\{\theta_j\}$,
$j\in\BZ/n\BZ$
are $\theta$-functions of order $n$.
Note that  here $\tau\in\BC$, we use the same symbol for
$\tau\in\CE$
and for some preimage of this point in $\BC$; $\theta_0,\ldots,
\theta_{n-1}$ can be considered as sections of a line bundle
of degree $n$ on $\CE$.
\par
The main properties of algebras
$Q_{n,k}(\CE,\tau)$.
\par
\noindent
(a)\enspace
$Q_{n,k}(\CE,\tau)$
is a graded algebra
$Q_{n,k}(\CE,\tau)=A_0\oplus A_1\oplus \ldots $~.
For generic $\tau$ is has the same size as the ring of polynomials.
It means that $\dim A_l=n(n+1)\ldots(n+l-1)/l!$.
If $\tau=0$ then
$Q_{n,k}(\CE,0)$
is an algebra of polynomials in $n$ variables.
So,
$Q_{n,k}(\CE,\tau)$
for small $\tau$ is a flat deformation of the space of functions
on $\BC^n$.
\par
\noindent
(b)\enspace
The finite Heisenberg group
$\Gamma_n$ (with generators
$\varepsilon_1,\varepsilon_2,\delta$,
$\varepsilon_1^n=\varepsilon_2^n=\delta^n=1$,
$\delta-$central element,
$\varepsilon_1\varepsilon_2=\delta \varepsilon_2\varepsilon_1$)
acts in
$Q_{n,k}(\CE,\tau)$
by automorphisms. This action is compatible with the grading,
and $A_1$ is an irreducible representation of
$\Gamma_n$.
\par
\noindent
(c)\enspace
For generic $\tau$ the center $Z$ of
$Q_{n,k}(\CE,\tau)$
is a polynomial ring  generated by $c$ elements
of degree $n/c$, where $c$ is the maximal common divisor
of $n$ and $k+1$.
\par
\noindent
(d)\enspace
{\it Characteristic manifold.\/}
It is the set of modules $M$ over
$Q_{n,k}(\CE,\tau)$
of the simplest possible type.
Such $M$ is graded, $M=M_0\oplus M_1\oplus \ldots$~,
$\dim M_i=1$,
and $M$ is generated by $M_0$.
Characteristic manifold has a structure of algebraic variety,
and we denote it by
$\Ch_{n,k}$.
To describe
$\Ch_{n,k}$
we need the decomposition of $n/k$ into continuous fraction.
\eq
&{n\over k}=n_1-{1\over n_2-{\textstyle 1\over \textstyle n_3-%
{~\atop \textstyle \ddots {\textstyle 1 \over \textstyle n_p}}}},
\quad
n_i\ge 2,
\quad
1\le i\le p.
\cr
\endeq
It is clear that such a decomposition is unique.
Now consider the  product of $p$ copies of $\CE$,
$\CE^{(p)}=\CE_1\times \CE_2\times\ldots\times \CE_p$.
Let $\xi$ be a line bundle on $\CE$ such that
$\deg\xi = 1$ (in this case $\dim H^0(\xi)=1$).
Then construct the line bundle
$\tilde \xi$ on $\CE^{(p)}$,
$\tilde \xi=\xi^{n_1+1}\boxtimes
\xi^{n_2+2}\boxtimes
\xi^{n_3+2}\boxtimes
\ldots \boxtimes
\xi^{n_{p-1}+2}\boxtimes
\xi^{n_{p}+1}$, $\boxtimes$
means exterior tensor product.
Denote by $\Delta_{i,i+1}$
the divisor on $\CE^{(p)}:
\Delta_{i,i+1}=\{(z_1,\ldots,z_p),
z_i+z_{i+1}=0\}$.
Let $\bar\xi$
be a bundle on $\CE^{(p)}$
which is equal to
$\tilde \xi-(\Delta_{1,2})-
(\Delta_{2,3})-\cdots-
(\Delta_{p-1,p})$.
It is easy to see that $\dim H^0(\bar\xi)=n$.
The bundle
$\bar\xi$
determines the map
$\CE^{(p)}\to\BC P^{n-1}$.
The image is isomorphic to the characteristic manifold.
\par
In this paper we examine these algebras and generalizations
using the vector bundles on elliptic curves.
To be more precise, the language of bundles on elliptic curves can
be used for description of symplectic leaves of such algebras.
Algebra
$Q_{n,k}(\CE,\tau)$
goes to abelian when $\tau\to 0$
so, the family of algebras
$Q_{n,k}(\CE,\tau)$
determines the hamiltonian structure on $\BC^n$.
We are interested in the symplectic leaves of this structure.
Let us formulate the result.
It is known that indecomposable bundles on an elliptic curve are labeled
by two integers $(n,k)$, $k>0$ and by a point of elliptic curve.
(See \S1).
For such a bundle
$\xi_{n,k}$,
$n$ is its degree and $k$ is its rank.
Let us fix an indecomposable bundle, and consider the moduli space
$\Mod(\xi_{0,1};\xi_{n,k})$
of $k+1$-dimensional bundles on elliptic curve with 1-dimensional subbundle
$(\nu ,\rho)$, $\dim\rho =1$, $\dim\nu =k+1$,
$\rho$ is trivial and $\nu/\rho\cong\xi_{n,k}$.
Suppose that
$n>0$ (otherwise the space
$\Mod(\xi_{0,1};\xi_{n,k})$
is empty).
It is easy to see that
$\Mod(\xi_{0,1};\xi_{n,k})\cong
P\big(\Ext^1(\xi_{n,k};\xi_{0,1})\big)$.
Here $P$ means projective space associated with vector space
$\Ext^1(\xi_{n,k};\xi_{0,1})$.
In other words,
$\Mod(\xi_{0,1};\xi_{n,k})$
is the space of exact sequences:
\eq
&0\to \xi_{0,1}\to Y \to \xi_{n,k}\to 0
\cr
\endeq
which are considered up to an isomorphism.
\par
The hamiltonian structure on $\BC^n$ arising in the classical limit
$\tau \to 0$ of
$Q_{n,k}(\CE,\tau)$
is homogeneous, so it determines the hamiltonian
structure on the projective space
$\BC\BP^{n-1}$.
Denote this structure by
$h_{n,k}(\CE)$.
\proclaim
Theorem 1.
The following decompositions of
$\BC\BP^{n-1}$
coincide:
(a) the decomposition into the union of symplectic leaves
of the structure
$h_{n,k}(\CE)$.
(b) $\BC\BP^{n-1}\cong \Ext^1(\xi_{n,k};\xi_{0,1}) \cong
\Mod(\xi_{0,1};\xi_{n,k})$.
This moduli space is the union of strata.
Each stratum corresponds to a type of
$k+1$-dimensional bundle.
\par
Hamiltonian structures on
$\BC\BP^{n-1}$
are the particular cases of the following general construction.
Let $G$ be a semisimple Lie group and $P$-a parabolic subgroup in
$G; \Ga, p$-Lie algebras of $G$ and $P$.
Algebra $\Ga$ is a sum
$n_-\oplus\Ga_0\oplus n_+$,
where $\Ga_0$ is a Levi subalgebra in $p$,
$n_-$ is the radical in $p$, $p=n_-\oplus \Ga_0$.
The standard bialgebra structure of $\Ga$, $\Ga \to \Lambda^2\Ga$
can be restricted on $p$.
It will be a non-trivial cocycle $p\to \Lambda ^2p$
which  gives us an elements $\beta$ from
$H^1(p,\Lambda^2p)$
or from
$\Ext^1_p(p^*, p)$.
Here we use the symbol
$p$ for the adjoint representation.
Consider now the moduli space
$\Mod(P,\CE)$
of $P$-bundles on $\CE$.
Let $\mu$ be a point of
$\Mod(P,\CE)$.
A tangent space
$T_\mu$ at $\mu$ is isomorphic to the space
$H^1\big(\CE,\mu(p)\big)$.
Here if $\rho$ is a representation of $p$
we denote by $\mu(\rho)$
the corresponding vector bundle on $\CE$.
The element $\beta\in H^1(p,\Lambda^2p)\cong
\Ext^1_p(p^*, p)$
determines the exact sequence of
$p$-modules $0\to p \to s\to p^*\to 0$
and the sequence of vector bundles on $\CE:
0\to\mu(p)\to \mu(s)\to\mu(p^*)\to 0$.
The boundary homomorphism in the exact cohomology
sequence is:
$H^0\big(\mu(p^*)\big)\to
H^1\big(\mu(p)\big)$.
Note that dual to
$H^0\big(\mu(p^*)\big)$
is
$H^1\big(\mu(p)\otimes K\big)\cong
H^1\big(\mu(p)\big)$,
because on elliptic curve the canonical bundle $K$
is trivial and we fix the trivialization.
Therefore we got a map
$\big(H^1(\mu(p))\big)^*\to
H^1\big(\mu(p)\big)$
which is skew-symmetric.
This map is determined at each point
$\mu\in \Mod(P,\CE)$,
so we construct a bivector field on
$\Mod(P,\CE)$
which is integrable.
Note that our construction  works on the non-singular part of
the manifold
$\Mod(P,\CE)$.
It is easy to see that the symplectic leaves of this structure are
exactly the fibers of the map
$\Mod(P,\CE)\to \Mod(G,\CE)$.
\par
Therefore it is a natural problem how to quantize such a hamiltonian
structure on the moduli space of $P$-bundles.
Algebras
$Q_{n,k}(\CE,\tau)$
give us a solution of this problem in a very special case.
In \S2 we consider the case when $P$ is a Borel subgroup in
a Kac-Moody group.
(Kac-Moody group can be infinite-dimensional).
(See point V of \S2).
The case when $P$ is a parabolic subgroup in
${\rm GL}_n$ will be the subject of our next paper.
\par
Now about the structure of the text.
In \S1 we collect  the facts about the  bundles
on elliptic curve.
The main problem which we deal with is: let $A$ and $B$ be two bundles
on $\CE$, then the space
$\Ext^1(B,A)$
has a natural stratification.
Each $\beta\in\Ext^1(B,A)$
determines the exact sequence
$0\to A\to \xi(\beta)\to B\to 0$.
The bundle $\xi(\beta)$ is a sum of indecomposable bundles.
We stratify the space $\Ext^1(B,A)$
according to the type of $\xi(\beta)$.
What is the combinatorial structure of this stratification?
In \S2 we present some constructions of algebras which in the simplest
case give us
$Q_{n,k}(\CE,\tau)$
when
$k=1$. In \S3 we give the description of symplectic leaves
in twisted version of
$Q_{n,k}(\CE,\tau)$
and discuss the case of a general parabolic group.

\beginsection
\S 1.
Indecomposable vector bundles on an elliptic curve

1.\ {\bf Description of the indecomposable bundles}

\noindent
Here we recall the results from [1].

\noindent
Construction of the stable indecomposable bundles.

Fix two integers $n,k, n\not =0, k>0$ which are relatively
prime and a line bundle $\alpha$ on $\CE$ of
degree $n\cdot k.$ If $n>0$ then
$\dim H^0(\al)=n\cdot k, H^1(\al)=0,$ if $n<0,$
then $H^0(\al)=0$ and $\dim H^1(\al)=n\cdot k.$
In both cases Heisenberg group
$\Gamma_{n\cdot k}$ (a central extension of
$\BZ_{nk}\times \BZ_{nk})$ acts on cohomologies in
irreducible way. The central element $K$ of $\Gamma_{n\cdot k}$
acts by multiplication on the primitive root of
unity
$\varepsilon, \varepsilon^{nk}=1$.
Consider Heisenberg subgroup
$\Gamma_k\hrw \Gamma_{nk}.$ The central element $K'$ of $\Gamma_k$
goes to $K^n.$ Consider also the irreducible representation
of $\Gamma_k, \Gamma_k\to {\rm End}\ \pi, $ where $K'$ acts by
$\varepsilon^{-k}$.
After tensoring $\al\otimes \pi$
we get a $k$-dimensional bundle on $\CE,$ where $\Gamma_k$
acts by the natural way (on $\al$ trough the
embedding $\Gamma_k\hrw \Gamma_{nk}$ and on $\pi$).
Central element acts by 1, so the group
$\BZ_k\times \BZ_k$ acts on $\al\otimes \pi.$
We can quote $\CE/\BZ_k\times\BZ_k$ and
$\al\otimes \pi/\BZ_k\times\BZ_k.$
As a result we get an indecomposable $k$-dimensional
vector bundle
$\xi_{n,k}(\al)$ on $\CE\cong\CE/\BZ_k\times\BZ_k.$
It is easy to see that
$\deg\xi_{n,k}(\al)=n, \dim H^0(\xi_{n,k}(\al))=n$
and $H^1(\xi_{n,k}(\al))=0$ if $n>0,$ and
$H^0(\xi_{n,k}(\al))=0$ and $\dim H^1(\xi_{n,k}(\al))=n$
if $n<0.$ The
Heisenberg group $\Gamma_n$ acts on cohomologies of
$\xi_{n,k}(\al)$ in irreducible way.
\hfill\break
\noindent
Construction of the semistable indecomposable bundles.

Fix again two integers $n,k$ such that
$n=n_1\cdot c, k=k_1\cdot c, k>0$
and $n_1$ and $k_1$ have not common divisors .
We denote by $\CO$ the trivial bundle on $\CE$ and
by $\CO(s)$ the bundle of $s$-jets of $\CO$.
It means that the fiber of $\CO(s)$ at the point
$a\in \CE$ is a quotient $\CO(a)/\CO(a)_s.$
Here $\CO(a)\cong \BC[[z]]$ is the space of sections of
$\CO$ in the formal vicinity of $a$ ($z$ is a local
coordinate in $a$) and $\CO(a)_s=z^s\BC[[z]]\subset \CO(a).$
It is easy to see that $\CO(s)$ is an indecomposable
$s$-dimensional bundle. We define $\xi_{n,k}(\al)$ as
$\xi_{n_1, k_1}(\al)\otimes\CO(c).$ If $n=0$ then by
$\xi_{0,k}(\al)$ we denote the bundle of $k$-jets of
a line bundle $\al, \deg \al=0.$

\noindent
This is the full list of indecomposable
bundles
$\{\xi_{n,k}(\al)\}; n,k\in \BZ$, $k>0$, $rk \xi_{n,k}(\al)=k$,
$\deg \xi_{n,k}(\al)=n$, $\al$ is a line bundle of degree
$n_1\cdot k_1$ where $n_1\cdot c=n, k_1\cdot c=k; n_1, k_1$
are relatively prime.
We will call the pair of integers $(n,k)$ the discrete
parameters of the bundle and $\al$-the continuous parameter.
In future we prefer to skip the continuous parameters and
write just $\xi_{n,k}$ if the dependence on $\al$ is evident.
In the next proposition we collect some simple properties of
$\xi_{n,k}$.

\proclaim
Proposition.

\noindent
(a) $\xi_{n,k}^*(\al)=\xi_{-n,k}(\al^{-1}).$
\hfill\break\noindent
(b) Suppose $n_2/k_2>n_1/k_1,$ then the dimension of
$\Hom(\xi_{n_1, k_1};\ \xi_{n_2, k_2})$ is
\hfill\break
$n_2\cdot k_1-n_1k_2;
\hfill\break
\Ext^1(\xi_{n_1, k_1};\ \xi_{n_2, k_2})=0$
\hfill\break\noindent
(c) $\Ext^i(\xi_{n_1, k_1};\ \xi_{n_2, k_2})
\cong \Ext^{1-i}(\xi_{n_2, k_2};\xi_{n_1, k_1})$
\hfill\break\noindent
(d) $\dim\Hom(\xi_{n,k}(\al_1); \xi_{n,k}(\al_2))=c, c$
is the maximal common divisor of $n$ and $k$ if
$\al_1=\al_2,$ otherwise it is zero.
\hfill\break
\smallskip
\noindent
2. {\bf The Moduli space of $k$-dimensional bundles on
an elliptic curve.}
\hfill\break\noindent
We know all indecomposable bundles, so each $k$-dimensional bundle
of degree $n$ can be decomposed into the sum:
\leq
&(*)\ \ \ \ \ \ \ \ \ \ \ \ \ \ \xi_{n_1,k_1}(\al_1)\oplus\xi_{n_2,k_2}(\al_2)
\oplus\ldots\oplus \xi_{n_p,k_p}(\al_p)\ \ \ \ \ \ \ \ \ \ \cr
\endeq
where
$
n_1+n_2+\ldots+n_p=n;
k_1+k_2+\ldots+ k_p=k.
$

So as a set the moduli space $\Mod_{n,k}$ of $k$-dimensional
bundles of degree $n$ is the set of data
$
\{ n_1, k_1, \al_1; n_2, k_2, \al_2; \ldots; n_p, k_p, \al_p\}.
$
\hfill\break
As usual we can organize the subset of semistable bundles
into some kind of a manifold.
(See for example [2])

\proclaim
Proposition.
A decomposition $(*)$ determines the semistable bundle if
and only if $n_1/k_1=n_2/k_2=\ldots=n_p/k_p.$
Such bundle is stable if and only if in each fraction
$n_i/k_i, n_i$ and $k_i$ have not common divisors.

We will use the symbol $\Mod_{n,k}^s$ for the moduli space of
stable bundles. It is clear that $\Mod_{n,k}^s$ is
isomorphic to the symmetric power $S^c\CE, c=(n,k).$
Recall that if we have an arbitrary family of vector
bundles parametrized by the manifold $N$ then the
semistable bundles constitute an open subset $U\subset N.$
Suppose that the dimension of bundles and the degree are
fixed. In this case we get the natural map $U\to \Mod_{n,k}^s.$
\hfill\break
\smallskip
\noindent
3. {\bf The Operator $F$}

$F$ is a well-known duality transform on the derived
category of the coherent sheaves on an elliptic curve.
Recall briefly the construction. Let $\CE'$ be the set of
line bundles on $\CE$ of degree zero, $\CE'\cong \CE.$
There is a natural bundle $P$ on the product $\CE\times \CE'.$
The fiber of $P$ at the point $(a,b)$ is the fiber of the bundle
$b$ at the point $a$. Let $\pi_1, \pi_2$ be two projections
$\CE\times\CE\to \CE:\pi_1(z_1, z_2)=z_1,
\pi_2(z_1,z_2)=z_2.$ Let $\nu$ be a coherent sheaf on $\CE.$
The transformation $\nu\to F(\nu)$ is given by the formula:
$F(\nu)=\pi_{2,*}\ (P\otimes\pi_1^*(\nu)).$
It is clear that $F$ is a covariant functor on the derived
category of coherent sheaves.
\hfill\break\noindent
{\bf Properties} of the functor $F.$
\hfill\break\noindent
(a) $F(\xi_{n,k}(\al))=\xi_{-k,n}(\al^{-1})$ if $n>0$
\hfill\break\noindent
(b) suppose that $n_1, n_2>0$, then
$\Hom(\xi_{n_1,k_1}; \xi_{n_2,k_2})
=\Hom(F(\xi_{n_1,k_1}); F(\xi_{n_2,k_2}))
\hfill\break\noindent
=\Hom(\xi_{-k_1,n_1}; \xi_{-k_2,n_2})
=\Hom(\xi_{k_2,n_2}; \xi_{k_1, n_1}).$
The same is true if we replace here $\Hom$ by ${\Ext}^1.$
\hfill\break
\noindent
4. {\bf Partial ordering on the moduli space of bundles.}
\hfill\break\noindent
Let $\xi_{n_1, k_1}(\al_1)\oplus \xi_{n_2, k_2}(\al_2)
\oplus\ldots\oplus\xi_{n_p,k_p}(\al_p)$
be a decomposition of some vector bundle $B$. We will say
that the bundle $C$ is less than $B$
\hfill\break\noindent
(notation: $C\prec B$)
if $C$ belongs to the closure of $B$. It means that
there is a
family $X$ of bundles which
satisfies the following properties.

The base $M$ of the family $X$ is connected.
Let $M_B$ and $M_\cc$ be the subsets in $M$ consisting
of points $p\in M$ such the the corresponding bundle
$X(p)$ is $B$ or $C$ respectively. Then $M_\cc$ belongs to
the closure of $M_B$.

It is clear
that the bundle $B$ is a maximal element (it means that if
$B\prec C$ then $C=B)$ if $\al_1, \al_2,\ldots , \al_p$
are all different and $n_1/k_1=n_2/k_2=\ldots=n_p/k_p$
and $n_i,k_i$ have not common divisors. Note, that in this
case the continuous parameters $\al_1, \ldots, \al_p$
are the line bundles of the same degree.
\hfill\break\noindent
{\it Problem\/}.
Suppose we have two bundles
$B=\xi_{n_1,k_1}(\al_1)\oplus\ldots\oplus\xi_{n_p,k_p}(\al_p)$
and
$
C=\xi_{\bn_1,\bk_1}(\beta_1)\oplus\ldots\oplus
\xi_{\bn_s,\bk_k}(\beta_s).
$
When $C\prec B$ ?

We restrict ourselves only by the case when all pairs
$(n_i,k_i)$
have not common divisors and the same is true for
the pairs $(\bn_i,\bk_i).$
Let us reformulate the question a little.
Let us fix a sequence of fractions
$\{\tau_1, \tau_2,\ldots,\tau_p\}, \tau_i=n_i/k_i, (n_i,k_i)=1,$
suppose that $\tau_1\le \tau_2\le\ldots\le \tau_p.$
We prefer to write $\xi_{\tau_i}(\al_i)$
instead of $\xi_{n_i,k_i}(\al_i).$
Suppose also that $\al_i\not=\al_j$ if $\tau_i=\tau_j.$
Then we deform the bundle
$\xi_{\tau_1}(\al_1)+\ldots+\xi_{\tau_p}(\al_p).$
The result is a $k_1+k_2+\ldots +k_p$ dimensional bundle $\nu$
which admit a filtration
$\nu_0\subset \nu_1\subset \ldots\subset \nu_{p-1}\cong \nu$
such that
$
\nu_0\cong \xi_{\tau_1}(\al_1); \
\nu_1/\nu_0\cong \xi_{\tau_2}\ (\al_2);\ldots ;
\nu_{p-1}/\nu_{p-2}\cong \xi_{\tau_p}\ (\al_p).
$
Informally, $p$ indecomposable bundles are glued together.
The existence of such a filtration is evident. For example,
after the deformation $\xi_{\tau_1}(\al_1)$ remains
to be a subbundle because
$\Ext^0(\xi_{\tau_i}(\al_i), \xi_{\tau_1}(\al_1))=0$
if, $i>1.$

The bundle $\nu$ can be decomposed into a sum of
indecomposable ones. Suppose, that all
indecomposable components of the decomposition are stable. So
$
\nu = \xi_{\bta_1}(\beta_1)\oplus\ldots
\oplus \xi_{\bta_s}(\beta_s),\bta_1\le
\bta_2\le\ldots\le \bta_s.
$

\noindent
What sequences $\{\bta_1,\bta_2,\ldots,\bta_s\}$
can be obtained this way?
To find an answer we need some preparations.
\hfill\break
\smallskip
\noindent
\maru{A}
Consider first the case of two bundles
$\xi_{\tau_1}(\al_1)$ and $\xi_{\tau_2}(\al_2),
0<\tau_1<\tau_2,$
$\tau_1=n_1/k_1; \tau_2=n_2/k_2.$
The space $W_{\tau_1,\tau_2}
=\Ext^1(\xi_{\tau_2}(\al_2), \xi_{\tau_1}(\al_1))$
has dimension $n_2k_1-k_2n_1.$
Each vector $\beta\in W_{\tau_1,\tau_2}$
determines the exact sequence:
\leq
&0\to \xi_{\tau_1}(\al_1)\to B[\beta]\to \xi_{\tau_2}(\al_2)\to 0.\cr
\endeq
Here $B[\beta]$ is a $k_1+k_2$-dimensional bundle which depends
on $\beta.$ So, we have a map
$W_{\tau_1,\tau_2}\ur{I} \Mod_{n_1+n_2, k_1+k_2}.$

Let us transform the exact sequence using the operator $F.$ We get:
\leq
&0\to \xi_{-\tau_1^{-1}}(\al_1^{-1})
\to FB[\beta]\to \xi_{-\tau_2^{-1}}(\al_2^{-1})\to 0,\cr
\endeq
In other words, there is a map:
$W_{-\tau_1^{-1},-\tau_2^{-1}}\ur{II}\Mod_{-k_1-k_2,n_1+n_2}$

\proclaim
Proposition.
The following diagram is commutative:
$$
\matrix{
W_{\tau_1,\tau_2}&\ur{I}&\Mod_{n_1+n_2, k_1+k_2}\cr
\apart
\da&&F\da\cr
\apart
W_{-\tau_1^{-1},-\tau_2^{-1}}&\ur{II}&\Mod_{-k_1-k_2,n_1+n_2}.\cr
}
$$

The fibers of maps I, II form  stratifications on the
linear spaces $W_{\tau_1,\tau_2}$ and
$W_{-\tau_1^{-1},-\tau_2^{-1}}.$
Our proposition claims that these stratifications are
identical.
Now fix some line bundle $\nu$ of degree $l.$
After tensoring $\xi_{\tau_1}(\al_1)$ and
$\xi_{\tau_2}(\al_2)$ by $\nu$ we get an isomorphism
$W_{\tau_1,\tau_2}\cong\Ext^1(\xi_{\tau_2}
(\al_2), \xi_{\tau_1}(\al_1))$ and
$W_{\tau_1+l, \tau_2+l} \cong \Ext^1 (\xi_{\tau_2+l}(\bal_2),
\xi_{\tau_1+l}(\bal_1)).$ Here $\bal_1$ and $\bal_2$
depend on $\al_1, \al_2$ and $\nu.$
It is clear that the fibers of the maps of $W_{\tau_1,\tau_2}$
and $W_{\tau_1+l, \tau_2+l}$ onto the corresponding
moduli spaces are the same. The result of these considerations
is the following. Let $\Gamma$ be the  group $SL_2(\BZ)$
which acts on the projective line in the usual way:
$
\pmatrix{
a&b\cr
c&d\cr
}
(\tau)\to {a\tau+b\over{c\tau+d}}, a,b,c,d\in \BZ,
ad-bc=1.
$
Then there is a commutative diagram:
$$
\matrix{
W_{\tau_1,\tau_2}&\rtw&\Mod,
\cr
\apart
\da&&\da\cr
\apart
W_{g(\tau_1), g(\tau_2)}&\rtw&\Mod\cr
\apart
}
$$
$g\in SL_2\ (\BZ).$
Therefore, the stratifications on the
spaces $W_{\tau_1, \tau_2}$ and $W_{g(\tau_1),g(\tau_2)}$
are the same.

We reformulate this is the next proposition.

\proclaim
Proposition.
Suppose that
$
\xi_{\tau_1}(\al_1)\oplus
\xi_{\tau_2}(\al_2)\prec \xi_{u_1}(\beta_1)\oplus
\xi_{u_2}(\beta_2)\oplus \ldots\oplus \xi_{u_s}(\beta_s).
$
Here $\tau_1, \tau_2, u_1,\ldots,u_s$ are rational numbers.
Then for each $g\in SL_2(\BZ)$
\leq
&\xi_{g(\tau_1)}(\bal_1)\oplus\xi_{g(\tau_2)}
(\bal_2)\prec\xi_{g(u_1)}(\bba_1)\oplus
\xi_{g(u_2)}(\bba_2)\oplus \ldots\oplus \xi_{g(u_s)}(\bba_s).\cr
\endeq
Note that continuous parameters change when we act by $g.$
\hfill\break
\smallskip
\noindent
\maru{B}
{\bf Rough stratification}

Start with a bundle
$
B=\xi_{\tau_1}(\al_1)\oplus
\xi_{\tau_2}(\al_2)\oplus\ldots\oplus
\xi_{\tau_p}(\al_p), \tau_1\le\tau_2\le\ldots\le
\tau_p,\tau_i=n_i/k_i, (n_i, k_i)=1.
$
Let us form a moduli space
$\Mod_{\tau_1, \tau_2,\ldots,\tau_p}$
which consists of the vector bundles of dimension
$k_1+k_2+k_3+\ldots+k_p$ with a filtration
$\nu_0\subset\nu_1\subset\nu_2\subset\ldots\subset\nu_{p-1}$
such that $\varphi_i:\nu_i/\nu_{i-1}\cong\xi_{\tau_{i+1}}
(\al_{i+1})$ and isomorphisms $\varphi_i$ are fixed.
The moduli space $\Mod_{\tau_1,\tau_2,\ldots,\tau_p}$
has a simple structure. Namely, it is a fibre bundle
$
\Mod_{\tau_1,\tau_2,\ldots,\tau_p}
\to \Ext^1(\xi_{\tau_2}(\al_2), \xi_{\tau_1}(\al_1))
\oplus\Ext^1(\xi_{\tau_3}(\al_3), \xi_{\tau_2}(\al_2))
\oplus\ldots\oplus\Ext^1(\xi_{\tau_p}(\al_p),
\xi_{\tau_{p-1}}(\al_{p-1}))).
$
Let ${\bM}$
be an arbitrary fiber of this bundle. Then it is
itself a fibre
bundle
$
\bM\to\Ext^1(\xi_{\tau_3}(\al_3), \xi_{\tau_1}(\al_1))
\oplus\Ext^1(\xi_{\tau_4}(\al_4), \xi_{\tau_2}(\al_2))
\oplus\ldots\oplus\Ext^1(\xi_{\tau_p}(\al_p),
\xi_{\tau_{p-2}}(\al_{p-2})).
$
Again we can consider the fiber of this bundle and it will
be a map onto the sum of finite dimensional vector spaces.
For the sequence of rational numbers
$\{s_1, s_2,\ldots, s_q\}, s_1\le s_2\le\ldots\le s_q$
we define a submanifold
$
\Mod
\matrix{
s_1,s_2,\ldots,s_q\cr
\tau_1,\tau_2,\ldots,\tau_p\cr
}
\subset \Mod_{\tau_1,\tau_2,\ldots,\tau_p}.$
A point $y\in \Mod_{\tau_1,\ldots,\tau_p}$
belongs to
$
\Mod
\matrix{
s_1,s_2,\ldots,s_q\cr
\tau_1,\tau_2,\ldots,\tau_p\cr
}
$
if the decomposition of the corresponding
$k_1+k_2+\ldots+k_p$ bundle is
$
\xi_{s_1}(\beta_1)\oplus
\xi_{s_2}(\beta_2)\oplus\ldots\oplus
\xi_{s_q}(\beta_q),
$
where $\beta_i\not=\beta_j$ if
$s_i=s_j.$ Manifold
$\Mod
\matrix{
s_1,\ldots,s_q\cr
\tau_1,\ldots,\tau_p\cr
}
$
is defined as the closure of the set of such $y.$
A group of symmetries of the bundle
$B$ is
$
\BC^*\times \BC^*\ldots\times\BC^*
$
($p$ times) acts on
$\Mod_{\tau_1,\tau_2},\ldots,\tau_p.$
The submanifolds
$
M
\matrix{
s_1,\ldots,s_q\cr
\tau_1,\ldots,\tau_p\cr
}
$
are stable with respect to this action. It is not
hard to show that
$M_{\tau_1,\ldots,\tau_p}$ is the union of all
$
M
\matrix{
s_1,\ldots,s_q\cr
\tau_1,\ldots,\tau_p\cr
}.
$
It is clear why we call this stratification ``rough''.
Suppose that
$s_1=s_2={p\over q}.$
Let us consider the manifold $U$ which is
$
M
\matrix{
s_1,\ldots,s_q\cr
\tau_1,\ldots,\tau_p\cr
}.
$
minus all strata which are smaller.
It is possible that $U$ contains the point such that
the corresponding vector bundle is
$\xi_{2p,2q}(\beta_1)\oplus\xi_{s_3}(\beta_3)
\oplus\ldots\oplus\xi_{s_q}(\beta_q).$
So we collect all such strata into one. The
bundles of ``generic'' type where all components in the
decomposition into a sum of indecomposable are stable
constitute the open set in the strata.
\hfill\break
\smallskip
\noindent
\maru{C}
{\it Characteristic manifold\/} and the ordering on the set
of sequences of rational numbers. The weak form of the
question about the ordering on the set of types of bundles
is following.

\proclaim
Definition.
Suppose $\{ \tau_1\le\tau_2\le\ldots\le\tau_p\}= R$
and
$\{s_1,\le\ldots,\le s_q\}= S$
are two sequences of rational numbers,
$\tau_i=n_i/k_i, s_j=\bn_j/\bk_j$ such that
$\Sigma n_i=\Sigma\bn_i$ and
$\Sigma k_i=\Sigma\bk_i.$
We will say that $R\prec S$ if the manifold
$
\Mod
\matrix{
s_1,\ldots,s_q\cr
\tau_1,\ldots,\tau_p\cr
}.
$
is non-empty.

\proclaim
Definition.
Characteristic manifold in
$\Mod_{\tau_1,\tau_2,\ldots,\tau_p}$
is the union of all strata
$
\Mod
\matrix{
s_1,s_2,\ldots,s_q\cr
\tau_1,\tau_2,\ldots,\tau_p\cr
}
$
such that if\
$
\Mod
\matrix{
s_1,s_2,\ldots,s_q\cr
\tau_1,\tau_2,\ldots,\tau_p\cr
}
\supset
\Mod
\matrix{
\bs_1,\bs_2,\ldots,\bs_l\cr
\tau_1,\ldots,\tau_p\cr
}
=Y
$
then $Y$ is a one-point stratum
\hfill\break\noindent
$
\Mod
\matrix{
\tau_1,\tau_2,\ldots,\tau_p\cr
\tau_1,\tau_2,\ldots,\tau_p\cr
}.
$
Therefore the characteristic manifold consists of the
``smallest'' possible strata. If
$R=\{ \tau_1,\ldots,\tau_p\},$ and
$S=\{s_1,\ldots,s_q\}$ is a stratum from the
characteristic manifold we will write down
$R\mapsto S.$

\proclaim
Proposition.
Sequence $R\prec S$ if and only if there is a set of sequences
$R\mapsto U_1\mapsto U_2\mapsto \ldots\mapsto U_l\mapsto S.$

\noindent
{\it Description\/} of the relation $\mapsto.$
Let $R=\{0,\tau\}$,
$\tau > 1$.
First we write down $\tau={n\over k}$
as a continuous fraction:
$n/k=(n_1,n_2,\ldots,n_t).$
The notation $(n_1,n_2,\ldots,n_t)$ means:
$$
(n_1,n_2,\ldots,n_t)
=n_1-{1\hfill\over n_2-{\textstyle 1\hfill\over \textstyle n_3-1%
{~\atop
\textstyle \ddots {~\atop{\textstyle -1 \over \textstyle n_t}}}}}
\qquad
\vtop{\hsize=5cm\parindent=0pt\overfullrule=0pt\strut
Here all $n_i \ge 2$.
The numbers $n_1,\ldots,n_t$ are determined by $n$, $k$.
\strut}
$$

\proclaim
Proposition.
The relation $\mapsto$ connects the pair
$\{0,\tau\}$ with the only one sequence
$\{ s_1, \ldots, s_t\},$ where
$s_1=1, s_2=n_1-1, s_3=(n_1, n_2-1),
s_4=(n_1, n_2,n_3-1),\ldots,
s_t=(n_1, n_2,\ldots,n_t-1).$

Now suppose that the pair $R$ is
$\{ \tau_1,\tau_2\},\tau_1\le \tau_2.$
There is again only one sequence
$S=\{s_1,\ldots,s_t\}$ such that
$R\mapsto S.$ To find it we can use the action of
$SL(2,\BZ).$ So there is $g\in SL(2,\BZ)$ such that
$g(\tau_1)=0$,
$g(\tau_2)>1$
and we first transform the pair
$\{\tau_1,\tau_2\}$ into $\{0,g(\tau_2)\},$
then go to the sequence $\{\bs_1, \bs_2,\ldots,\bs_t\}$
according to the previous proposition and after
that we return back with help of
$
g^{-1}: \{s_1,\ldots,s_t\}=\{g^{_1}(\bs_1),
g^{-1}(\bs_2),\ldots,g^{-1}(\bs_t)\}.
$
Combining all this arguments we get the next proposition.

\proclaim
Proposition.
The relation $\mapsto$ connects the pair
$\{\tau_1,\tau_2\}$ with the only one sequence
$\{s_1,s_2,\ldots,s_t\}$ which is calculated in the
following way. Suppose first that the integegral part of
$\tau_1$ is less then the integegral part of $\tau_2.$
Write down
$\tau_1$ as $m_1+(m_2,m_3,\ldots,m_a)^{-1}$ and
$\tau_2$ as $(n_1,n_2,\ldots,n_b).$ Then
$\{s_1,s_2,\ldots,s_t\}$ is
$
\{m_1+(m_2,m_3,\ldots,m_a-1)^{-1},
m_1+(m_2,m_3,\ldots,m_{a-1}-1)^{-1},
\ldots,m_1+(m_2-1)^{-1}, m_1+1,
n_1-1, (n_1, n_2-1), (n_1,n_2,n_3-1),
\ldots, (n_1,n_2,\ldots,n_b-1)\}.
$
Here we suppose that $m_2,m_3,\ldots \ge 2.$

General case is clear because of the acting
$SL_2(\BZ).$

\remark
Suppose that $\tau_1=p_1/q_1, \tau_2=p_2/q_2$
such that $p_2q_1-q_2p_1=1$ then the sequence
$\{s_1\ldots s_t\}$ consists only of one term. So
$\{\tau_1,\tau_2\}\mapsto\{(p_1+p_2)/(q_1+q_2)\}.$
Now we describe relation $\mapsto$ in the general case.

\proclaim
Proposition.
The sequences
$\{\tau_1,\ldots,\tau_p\}\mapsto
\{s_1,s_2,\ldots ,s_q\}$ if the sequence
$\{s_i\}$ can be obtained by the following construction.
Let us fix some pair of $\{\tau_i, \tau_{i+1}\}.$
Then using the rule from the previous proposition
write the sequence
$S=\{f_1,\ldots,f_t\}, \{\tau_i,\tau_{i+1}\}\mapsto S.$
\hfill\break\noindent
Then the sequence $\{s_1,s_2,\ldots,s_q\}$ will be
\leq
&\{\tau_1, \tau_2,\ldots, \tau_{i-1}, f_1, f_2,\ldots,
f_t, \tau_{i+2},\ldots, \tau_p\}.\cr
\endeq

Using this rule it is possible to find all rough strata
in the manifold $\Mod_{\tau_1,\tau_2,\ldots,\tau_s}.$
They are labeled by the sequences
$\{s_1,\ldots,s_q\}$ which are smaller than
$\{\tau_1,\ldots,\tau_s\}.$
All such sequences can be discribed with help of the
system of inequalities. We will write them in the special
case of $(0,\tau)$-pair.

\proclaim
Proposition.
Let $\tau$ be $n/k$, $(n,k)=1$, then the sequences
$S=\{p_1/q_1, p_2/q_2,\ldots, p_l/q_l\}$,
$(p_i,q_i)=1$, $p_i/q_i\le p_{i+1}/q_{i+1}$ for
$i=1,2,\ldots,l-1$ which satisfy the
condition $(0,\tau)\prec S$ can be characterized as a
solution of the following system:
\hfill\break
\noindent
$
0<p_1/q_1, p_l/q_l<n/k, p_1>0,\
p_2+p_1q_2-p_2q_1>0,\
p_3+(p_1+p_2)q_3-p_3(q_1+q_2)>0\ldots
p_l+(p_1+p_2+\ldots+p_{l-1})
q_l-p_l(q_1+q_2+\ldots+q_{l-1})>0;
$
\hfill\break
$
p_1+\cdots+p_l=n,\
q_1+\cdots+q_l=k+1$.

\proclaim
Example (1).
The rough strata which are more than the fixed
one constitute the set with partial ordering.

\item{(a)}
$N$ is even, then
$\{0,N\}\mapsto\{1,N-1\}
\mapsto\ldots\mapsto\{{N\over 2}, {N\over 2}\}$

\item{(b)}
$N$ is odd then
$\{0,N\}\mapsto\{1,N-1\}
\mapsto\ldots
\{{{N-1}\over 2}, {{N+1}\over 2}\}
\mapsto\{{N\over 2}\}$

\item{(c)}
$\{0, {5\over2}\}\mapsto\{1,2,2\}\mapsto
\{{3\over2},2\}\mapsto\{{5\over3}\}$

\item{(d)}
$
\{0, {17\over2}\}\mapsto\{1,8,8\}\mapsto\{2,7,8\}\yaya
$
$$
\matrix{
\mapsto&\{3,6,8\}&\mapsto&\{4,5,8\}&\mapsto&
\{{9\over 2},8\}&\mapsto&\{5,5,7\}&
\mapsto&\{5,6,6\}&\to\cr
\apart
&&\searrow&&\searrow&&\nearrow&&\nearrow\cr
\apart
\mapsto&\{2,{15\over2}\}&\mapsto&
\{3,7,7\}&\mapsto&\{4,6,7\}&\mapsto&\{4,{13\over3}\}\cr
}
$$
$
\matrix{
\mapsto&\{{11\over 2},6\}&\mapsto&\{{17\over3}\}.\cr
}
$

\item{(e)}
In this point we start with a pair $\{0,(n_1,n_2)\}$
where $n_1,n_2$ are sufficiently big.
\hfill\break
\noindent
$
\{0,(n_1,n_2)\}\mapsto\{1, n_1-1, (n_1, n_2-1)\}
\yaya \hbox{(to the next line)}
$
$$
\matrix{
\mapsto&\{2, n_1-2, (n_1,n_2-1)\}&\mapsto&
\{3, n_1-3, (n_1,n_2,-1)\}&\ayaya\cr
&&\searrow\cr
\mapsto&\{1, n_1-{1\over2},(n_1,n_2-2)\}&
\mapsto&\{2, n_1-1,n_1-1, (n_1,n_2-2)\}&\ayaya\cr
&&\searrow\cr
&&&\{1, n_1-{1/3}, (n_1,n_2-3)\}&\ayaya\cr
}
$$

$$
\matrix{
&&&\{5,n_1-5,(n_1,n_2-1)\}\hfill\cr
&&\nearrow\cr
&\{4,n_1-4,(n_1,n_2-1)\}\hfill&\mapsto&\{4,n_1-3,n_1-1,(n_1,n_2-2)\}\hfill\cr
\ayaya&&\nearrow\cr
&\{3,n_1-2,n_1-1,(n_1,n_2-2)\}\hfill&\mapsto&\{3,n_1-{3\over
2},(n_1,n_2-2)\}\hfill\cr
\nearrow&&\searrow\cr
&&&\{3,n_1-2,n_1-{1\over 2},n_1-{1\over n_2-3}\}\hfill\cr
\searrow&&\nearrow\cr
&\{2,n_1-1,n_1-{1\over
2},(n_1,n_2-3)\}\hfill&\mapsto&\{2,(n_1,2,2),(n_1,n_2-3)\}\hfill\cr
\ayaya&&\searrow\cr
&\{1,n_1-{1\over 4},(n_1,n_2-4)\}\hfill&\mapsto&\{2,n_1-1,n_1-{1\over
3},(n_1,n_2-4)\}\hfill\cr
&&\searrow\cr
&&&\{1,n_1-{1\over 5},(n_1,n_2-5)\}\hfill\cr}
$$

\proclaim
Example (2).
Here we consider $(9/2)$ case. The space
$\Ext^1(\xi_0, \xi_{9/2})$
is decomposed according to the type of
corresponding bundle. So $\BC^9$ is a union of
submanifolds.

\item{(a)}
origin where decomposition is $\xi_0\oplus\xi_{9/2}$

\item{(b)}
characteristic manifold. It has dimension three,
and the bundle is a sum
$\xi_1(\al_1)\oplus\xi_4(\al_2)\oplus\xi_4(\al_3)$,
if $\al_3\not=\al_4$ and $\xi_1(\al_1)\oplus\xi_{4/2}(\al_2)$
if $\al_2=\al_3$,
where the sum $\al_1+\al_2+\al_3=\mu$ is fixed and
equal to the determinant of the bundle $\xi_{9/2}$.

\item{(c)}
next piece is a 5-dimensional manifold, where decomposition
is $\xi_2(\al_1)\oplus\xi_3(\al_2)\oplus\xi_4(\al_3),
\al_1+\al_2+\al_3=\mu.$ The set which corresponds to
fixed $\al_1,\al_2,\al_3$ has dimension 3.

\item{(d)}
two 6-dimensional manifolds where decomposition is
$\xi_2(\al_1)\oplus\xi_{7/2}(\al_2), \al_1+2\al_2=\mu$ and
$\xi_{5/2}(\al_1)\oplus\xi_4(\al_2), 2\al_1+\al_2=\mu$.

\item{(e)}
$9$-dimensional part, where decomposition is
$\xi_3(\al_1)\oplus\xi_3(\al_2)\oplus\xi_3(\al_3)$,
$\al_1+\al_2+\al_3=\mu.$
If $\al_1,\al_2,\al_3$ are all different then
corresponding submanifold in $\BC^9$ has dimension 9.
Suppose $\al_1=\al_2.$ Then for fixed
$\al_1=\al_2, \al_3\not=\al_1$ the submanifold has
dimension 5. The stratum where we have
$\xi_{6/2}(\beta)\oplus\xi_3(\al)$
has dimension 6.
Here $\beta$ and $\al$ are chosen such that
 $\Hom(\xi_{6/2}(\beta), \xi_3(\al))=0.$
If $\Hom(\xi_{6/2}(\beta), \xi_3(\al))=\BC$
then the dimension is 4. The last possibility
$\xi_3(\al)\oplus\xi_3(\al)\oplus \xi_3(\al).$ Here
$3\al=\mu$ and for each such $\al$
(we have $9$ possibilities) it is a line in $\BC^9.$

\beginsection
\S 2. Construction of the algebras

I. The ``symmetric'' algebra associated with $R$-matrix (due to Lusztig
$\cdots$)

Fix the vector space $V$ and $R$-matrix $V\otimes V\to V\otimes V$.
We define first the so-called ``exchange'' algebra associated with
$R$. Consider the algebra $A$ with the space of generators $V[i],i\in \BZ$.
Each $V[i]$ is isomorphic to $V$. Note, that the $R$-matrix defines a subspace
$W\in V\otimes V\oplus V\otimes V$ which consists of pairs $v_1\otimes
v_2-R(v_1\otimes v_2)$;
let $W_{ij}$ be the same space in $V[i]\otimes V[j]\oplus
V[j]\otimes V[i]$. The defining relations in $A$ are: (a) $V[i]\otimes V[i]$,
it means that $v_1v_2=0$ if $v_1,v_2\in V[i]$. (b) $W_{ij}\subset V[i]
\otimes
V[j]\oplus V[j]
\otimes
V[i]$, $i<j$. It is easy to see that as a vector space
$A=\oplus V[i_1]\otimes V[i_2]\otimes \cdots \otimes V[i_n]$
where the sum goes over all sequences $i_1<i_2<\cdots <i_n$. If we
replace the direct sum by the product we get an algebra
$\bar A=\sqcap V[i_1]\otimes \cdots \otimes V[i_n]$.
The element from $\bar A$ is a function from the set of sequences
$\{ i_1<i_2<\cdots <i_n\}$ into $V\otimes V\otimes \cdots \otimes V$
(n times). The product of such infinite expressions is well-defined.
Now let us construct the maps $\theta_n:
\underbrace{V\otimes \cdots \otimes V}
_{n~{\rm times}}
\to \bar A$.
For $v\in V$ we denote by $v[i]$ the corresponding element from $V[i]$.
The maps are: $\theta_n(v_1\otimes v_2\otimes \cdots \otimes v_n)
=\sum v_1[i_1]\cdot v_2[i_2]\cdot \cdots \cdot v_n[i_n], i_1<i_2<\cdots <i_n$.
Denote by $T(V)$ the direct sum $\BC\oplus V\oplus V\otimes V\oplus \cdots,
\{\theta_n\}$ constitute a map $T(V)\to \bar A$ and the image is a subalgebra.
Our notation for $T(V)$ with this multiplicative structure will be
$T_R(V)$. A ``Symmetric'' algebra associated with $R$-matrix is a
sub-algebra $S_R(V)\subset T_R(V)$ which is generated by $V\subset T(V)$.
If $R$ is just permutation the $T_R(V)$ is a commutative algebra and $S_R(V)$
is a polynomial ring. If $R$ is unitary in the sense that $R^2=1$
then $S_R(V)$ has the same size as a polynomial ring. In the case when
$R^2\not=1$,
$S_R(V)$ is bigger than a polynomial algebra, and for ``generic''
$R$, $S_R(V)$ is isomorphic to a free algebra.

Algebras $T_R(V)$ and
$S_R(V)$ bear the natural structure of the braided Hopf algebras.
Namely fix in $\bar A$ two subalgebras $T_1$ and $T_2$. The first one
is spanned by
the elements $\sum v_1[i_1]\cdot v_2[i_2]\cdot \cdots \cdot v_n[i_n],
i_1<i_2<\cdots <i_n\le 0$. The second has a basis $\sum v_1[i_1]\cdot v_2[i_2]
\cdot \cdots \cdot v_m[i_m], 0<i_1<i_2<\cdots <i_m$. It is clear that
$T_1$ and $T_2$ are isomorphic to $T_R(V), T_1$ and $T_2$ generate
in $\bar A$ the subalgebra $T_R(V)\bar \otimes T_R(V)$. This algebra
has the same size as a usual tensor product $T_R(V)\otimes T_R(V)$,
but the multiplication is twisted in the evident way. (For example,
the spaces $V\otimes 1$ and $1\otimes V$ from the $T_R(V)\bar \otimes T_R(V)$
are not commute but satisfy $R$-matrix permutation relation).
$T_R(V)\bar \otimes T_R(V)$ is called the braided tensor product
of two copies of $T_R(V)$. Algebra $T_R(V)$ is a subalgebra in
$T_1\bar \otimes T_2$, therefore we get a comultiplication map
$T_R(V)\to T_R(V)\bar \otimes T_R(V)$. The same arguments can be applied to
$S_R(V)$. As a result $S_R(V)$ has a structure of braided Hopf algebra.
It means that
there is a comultiplication map $\Delta :S_R(V)\to S_R(V)\bar \otimes
S_R(V)$ where $\bar \otimes$
is the twisted tensor product and $\Delta$ on the
generators is given by the formula $\Delta(v)=v\otimes 1+1\otimes v$.

\proclaim
Example 1. Consider the case when $R$ is diagonal. So $V$ has a basis
$\{ v_1,v_2,\cdots ,v_n\}$ such that $R(v_i\otimes v_j)=C_{ij}(v_j\otimes
v_i), C_{ij}\not= 0$. If matrix $\{ C_{ij}\}$ satisfies the symmetry
condition $C_{ij}C_{ji}=1$ then the algebra $S_R(V)$ is an algebra of skew
polynomials with generators $V_i$ and relations $v_iv_j=C_{ij}v_jv_i$. If
$C_{ij}=B_{ij}\cdot D_{ij}$ where $D_{ij}D_{ji}=1$ then the ``symmetric''
algebras which correspond to $\{ C_{ij}\}$ and $\{ B_{ij}\}$ are very
close. More precisely: let $S(C)$ will be the first and $S(B)$-the second.
All algebras $S(C), S(B)$ and
the algebra of skew polynomials $S(D)$ have a
$\Gamma$-grading, where $\Gamma$ is a lattice with basis $\{ e_i\}; \deg v_i
=e_i$. Then $S(C)$ is a subalgebra in $S(B)\otimes S(D)$ which is
generated by the products $v_i\otimes v'_i$
($v_i$-generators of $S(B)$, $v'_i$-of
$S(D))$. It is evident the $S(C)$ has the same size as $S(B)$ and
coincides with
the subalgebra in $S(B)\otimes S(D)$ of elements of degree zero
($\deg v_i=e_i,\deg v'_i=-e_i$). This all means that the most interesting case
is when the matrix $\{ C_{ij}\}$ is symmetric. In this case if we write down
$C_{ij}=q^{c_{ij}}$, where $\{c_{ij}\}$ is a matrix of scalar products
of the system of simple roots of some Kac-Moody algebra,
then algebra $S_R(V)$
is isomorphic to the universal enveloping algebra of the maximal nilpotent
subalgebra in the $q$-deformed Kac-Moody algebra.

\bigskip
\noindent
II. Algebra which is defined by diagonal $R$-matrix in ``continuous'' case.

We apply now the scheme of the point I in the case when $R$ is diagonal but
$V$ is the  space of functions in one variable. If $V$ is the space of
polynomials then $V\otimes V$ is the space of polynomials in two variables.
We do not want to specify what functions are in $V$. And we denote by
$V\tilde \otimes V$ the space of functions in two variables. Now fix a
function $\lambda(x,y)$ in two variables and
let $R$-matrix $R_\lambda$
acts $V\tilde \otimes V\to V\tilde \otimes V$ by the formula
$R_\lambda(g(x,y))=g(y,x)\cdot \lambda(x,y)\cdot \lambda^{-1}(y,x),
g(x,y)\in V\tilde \otimes V$. Formally $R_\lambda\cdot R^{-1}_\lambda=1$,
but the algebra $S_{R_\lambda}(V)$ can be more complicated than the ring
of skew polynomials. The reason is that we have to replace the usual tensor
product by the ``extended'' one $\tilde \otimes$. First we construct the
bigger algebra $\tilde S_{R_\lambda}(V)=\BC\oplus S_1\oplus S_2\oplus \cdots ,
S_j$ consists of the symmetric functions in $j$ variables. The product
is given  by the formula: $f\in S_\alpha , g\in S_\beta$
\leq
&f*g(X_1,\cdots ,X_{\alpha+\beta})\cr
&={1\over {\alpha !\cdot \beta !}}
Symm \Bigg[\bigg[
\prod_{{1\le i\le \alpha}\atop{\alpha <j\le \alpha+\beta}}
\lambda(X_i,X_j)\bigg]
\cdot f(X_1,\cdots ,X_\alpha)\cdot g(X_{\alpha+1}\cdots ,
X_{\alpha+\beta})\Bigg].\cr
\endeq
Here symbol Symm means the symmetrization with respect to
the symmetric
group acting on the set of $\alpha+\beta$ variables.

\proclaim
Definition. Algebra $S_{R_\lambda}(V)$ (to simplify the notations $S_\lambda$)
is the subalgebra in $\tilde S_{R_\lambda}(V)$ generated by $S_1$.

\proclaim
Examples. If the function $\lambda={{x-qy}\over{x-y}}$ then the algebra
$S_\lambda$
coincides with a subalgebra in the universal enveloping of $q$-deformed
$\hat{sl}_2$. This subalgebra is a $q$-version of
the universal enveloping algebra of the
algebra of currents on the line into the maximal nilpotent subalgebra of
$sl_2$.

We need also the shifted version of the ``symmetric'' algebra. Suppose the
$R$-matrix $R:V\otimes V\to V\otimes V$ admit symmetries. In the simplest
case it means that some operator $A: V\otimes V$ has a property
$R\circ(A\otimes A)=(A\otimes A)\circ R$. Operator $A$ is extended to
the automorphism
$A:S_R(V)\to S_R(V)$. In such a situation we can add an element
$a$ to the algebra $S_R(V)$ with the permutation relations:
$a\cdot u=A(u)\cdot a$, where $u\in S_R(V)$. As a result we get semi-direct
product $\BC [a]\ltimes S_R(V)$. The shifted algebra $S_{R,A}(V)$ is by
definition a subalgebra in $\BC[a]\ltimes S_R(V)$ which is generated by
the subspace $a\cdot V$. It is clear that $S_{R,A}(V)$ has the same size as
$S_R(V)$. Let us apply this construction in the case, when $V$ is the
space of functions on the line and $A$ is a shift $f(z)\to f(z-p)$.
This operator commutes with diagonal $R$-matrix if $\lambda(x,y)=
\lambda(x+p,y+p)$. Suppose $\lambda$ satisfies this condition, then
the product in $S_\lambda$ is deformed in
the following way:
\leq
&f*g(x_1,x_2,\cdots ,x_{\alpha +\beta})\cr
&={1\over{\alpha !\beta!}}
Symm \Bigg(\bigg(
\prod_{{1\le i\le \alpha}\atop{\alpha <j\le \alpha + \beta}}
\lambda(x_i,x_j)\bigg)
f(x_1,\cdots ,x_\alpha)g(x_{\alpha+1}-\alpha p,x_{\alpha+2}
-\alpha p,\cdots ,x_{\alpha+\beta}-\alpha p)\Bigg).\cr
\endeq
We denote the ``shifted'' algebra
by $S_{\lambda,p}$.

\bigskip
\noindent
III.
\proclaim
Proposition. Suppose that $\tau$ is not a point of finite order in $\cal E$.
The algebra $Q_{n,1}(\CE ,\tau)$ is a subalgebra in $S_{\lambda,p}$,
where $p=2\tau$ and $\lambda(x,y)
=\theta(x-y-n\tau)/\theta(x-y)$. The space of
generators of $Q_{n,1}(\CE,\tau)$ is the space of $\theta$-functions
of order $n$.

Here $\theta$- is a $\theta$-function of order one which
has zero at the
origin. The meaning of this proposition is that the algebra
$Q_{n,1}(\CE,\tau)=
\BC\oplus B_1\oplus B_2\oplus \cdots$ has the following realization.
The space of generators $B_1$ is identified with the space
$H^0(\xi)$, where $\xi$ is a line bundle of degree $n$. The space $B_j$ is
isomorphic to the space of sections of the line bundle
$S^p\xi^{(1-p)\tau}$. Here $\xi^{(1-p)\tau}$ is the line bundle $\xi$
shifted by the map $\CE \to \CE, x\to x+(p-1)\tau;
S^p\xi^{(1-p)\tau}$ is a symmetric power of $\xi^{(1-p)\tau}$, therefore
the sections of $S^p\xi^{(1-p)\tau}$ are symmetric $\theta$-functions of
order $n$ in $p$ variables.

Now we present the construction of the algebra $Q_{n,1}(\CE,\tau)$
which is the version of the constructions from point I (with some
modifications).

Let $A$ be an algebra with two generators $e,u$ and relations $e\cdot u
=(u+(n-2)\tau)\cdot e, \tau\in \BC$. Define the braided tensor product of
two copies of $A: A_1$ and $A_2$, namely, $A_1\bar \otimes A_2$ is an
algebra with generators $e_1,u_1;e_2,u_2$ and relations:
$e_1u_1=(u_1+(n-2)\tau)e_1, e_2u_2=(u_2+(n-2)\tau)e_2, e_1u_2=(u_2-2\tau)e_1,
e_2u_1=(u_1-2\tau)e_2, [u_1,u_2]=0,$ and $e_1 e_2=
-\exp(2\pi i(u_2-u_1))
{{\theta(u_1-u_2-n\tau)}\over {\theta(u_2-u_1-n\tau)}}\cdot e_2 e_1.$
We have two maps $\mu_1:A\otimes A\cong (A_1\otimes 1)\cdot (1\otimes A_2)
\to A_1\bar \otimes A_2$ and $\mu_2: (1\otimes A_2)\cdot (A_1\otimes 1)
\cong A\otimes A\to A_1\bar \otimes A_2$. The map $\mu_1\circ \mu^{-1}_2:
A\otimes A\to A\otimes A$ is a sort of $R$-matrix and it is possible to
apply the technique from I. The result is that $Q_{n,1}(\CE,\tau)$
is a subalgebra in $A\bar \otimes A\bar \otimes \cdots \bar \otimes A$.

\proclaim
Proposition. Let $B$ be an algebra with generators
$\{ e_\alpha,u_\alpha \}, \alpha \in I$,
(the set of indices I can be infinite). The relations are:
$e_\alpha u_\beta=(u_\beta-2\tau)e_\alpha, (\alpha \not= \beta),
e_\alpha u_\alpha
=(u_\alpha+(n-2)\tau)e_\alpha, e_\alpha e_\beta=-e^{-2\pi i(u_\beta-u_\alpha)}
{{\theta(u_\alpha-u_\beta-n\tau)}\over{\theta(u_\beta-u_\alpha-n\tau)}}
e_\beta e_\alpha, u_\alpha u_\beta=u_\beta u_\alpha$.
There is a homomorphism from the algebra
$Q_{n,1}(\CE,\tau)$ into $B$.
The generator $f$ from
$Q_{n,1}(\CE, \tau)$ (which is $\theta$-function of order $n$) goes to
$\sum\limits_\alpha
f(u_\alpha)\cdot e_\alpha \in B$.

The realization of $Q_{n,1}(\CE,\tau)$ as a space of symmetric
$\theta$-functions we call the {\it functional realization.} The center
of $Q_{n,1}(\CE,\tau)$ has a nice description in these terms.

\proclaim
Proposition.
(a) Suppose $n$ is even. Then the center $Z$ of algebra
$Q_{n,1}(\CE,\tau)$ is a polynomial algebra with two generators
of degree $n/2$. This space of generators coinsides with the subspace
$W_1\subset S^{n/2}\xi^{(1-n/2)\tau}$ which consists of $\theta$-functions
$g(z_1,\cdots ,z_{n/2})$
such that
$g(z_1,z_1-n\tau,z_3,\cdots ,z_{n/2})=0.$ The
quadratic expressions of generators constitute a 3-dimensional subspace
$W_2\subset S^n\xi^{(1-n)\tau};W_2$ is a space of functions
$g(z_1,\cdots ,z_n)$ with condition $g(z_1,z_1-n\tau, z_1-2n\tau,
z_4,\cdots ,z_n)=0; W_s\subset S^{ns/2}\xi^{(1-ns/2)\tau}$ is
a space of $g(z_1,\cdots ,z_{ns/2}),$ the condition is:
$g(z_1,z_1-n\tau,
\cdots ,z_1-sn\tau,z_{s+2}, \cdots ,z_{ns/2})=0$.
\hfil\break
(b) Suppose $n$ is odd. Then the center $Z$ of the algebra $Q_{n,1}(\CE,\tau)$
is generated by one element $\Delta$ of degree $n$, $\Delta$ is represented
by the function $g(z_1,\cdots ,z_n)$ which is zero on the diagonal $(z_1,
z_1-n\tau, z_1-2n\tau,z_4,\cdots ,z_n)$. There is only one function (up
to a constant) with this property. Note that element $\Delta^l$ is
a function in $n\cdot l$ variables $z_1,\cdots ,z_{nl}$. It is zero
if $z_1=z_2+n\tau, z_2=z_3+n\tau, \cdots z_{2l}=z_{2l+1}+n\tau$.

\bigskip
\noindent
IV. Generalization.
\par
We apply now our construction ``with function $\lambda$''
to find new algebras. Let $\Ga$ be the following Lie superalgebra:
$\Ga=\Ga_1\oplus \Ga_2, \Ga_1$ is odd part, $\Ga_2$ is even, $[\alpha,\beta]
=0$ if $\alpha,\beta\in \Ga_2$ or $\alpha\in\Ga_1, \beta\in \Ga_2$.
We identify
the space $\Ga_1$
with the space of sections of the line bundle
$\xi,\deg \xi=n$, and $\Ga_2=H^0(\xi^2)$. So, $\dim \Ga_1=n, \dim \Ga_2
=2n$. The bracket $
\Ga_1\otimes \Ga_1\to \Ga_2$ is the product of the sections.
The universal
enveloping algebra
of $\Ga$ is a graded algebra ($\deg \Ga_1=1, \deg \Ga_2=2$) and
its Hilbert function is equal to $(1+t)^n(1-t^2)^{-2n}=(1-t)^{-n}
(1-t^2)^{-n}$.
We want to construct now the elliptic deformation $T_n(\CE,\tau)$
of the $U(\Ga)$.
The algebra $T_n(\CE,\tau)=\BC\oplus A_1\oplus A_2
\oplus \cdots$, where $A_1=H^0(\xi)$, $A_p$ is the space of meromorphic
sections
of the bundle $\Lambda^p\xi^{-(p-1)\cdot \tau}$ on $S^p\CE
=\CE^p/S_p$ which satisfy the properties: 1) $f(z_1,\cdots ,z_p)$
is skew-symmetric and
has a pole of order $\le 1$ on the diagonal $z_1=z_2$ and is
holomorphic outside
the diagonal; 2) $f(z_1,\cdots ,z_p)=0$ if $z_1=z_2+n\tau=z_3+2n\tau$.
The multiplication of two elements $f\in A_\alpha, g\in A_\beta$ is
given by the formula:
\leq
&f*g(x_1,x_2,\cdots ,x_{\alpha+\beta})\cr
&={1\over{\alpha !\beta !}}Symm
\bigg(
\prod_{{1\le i\le\alpha}\atop{\alpha<j\le \alpha+\beta}}
{{\theta(x_i-x_j+2n\tau)\theta(x_i-x_j-n\tau)}\over{\theta(x_i-x_j)^2}}
\bigg)
f(x_1,\cdots ,x_\alpha)\cr
&\cdot g(x_{\alpha+1}+2\alpha \tau,x_{\alpha+2}
+2\alpha\tau,\cdots ,x_{\alpha+\beta}+2\alpha\tau)\cr
\endeq
Here $Symm$ is the operator of skew-symmetrization.
It is a homomorphism from $T_n(\CE,\tau)$ into the algebra with
generators $\{ u_\alpha,e_\alpha\}, \alpha\in I$
($I$-an arbitrary set of indices)
and relations: $[u_\alpha,u_\beta]=0; e_\alpha u_\beta=(u_\beta
+2\tau)e_\alpha$, in this relation $\alpha\not=\beta; e_\alpha u_\alpha
=(u_\alpha-(2n-2)\tau)e_\alpha;$
\leq
&e_\alpha e_\beta=-e^{4\pi i(u_\beta-u_\alpha)}\cdot
{{\theta(u_\alpha-u_\beta+2n\tau)\theta(u_\alpha-u_\beta-n\tau)}\over
{\theta(u_\beta-u_\alpha+2n\tau)\theta(u_\beta-u_\alpha-n\tau)}}
e_\beta e_\alpha\cr
\endeq
The generator $f\in A_1=H^0(\xi)$ goes to the sum $\sum_\alpha f(u_\alpha)
\cdot e_\alpha$.

\noindent
{\it Remark 1.} Our definition of the algebra $T_n(\CE,\tau)$
works in the
case when $\tau$ is not a point of finite order in $\CE$. If $\tau$
is a point of finite order we can define $T_n(\CE,\tau)$ as a
limit. But the procedure of finding out the limit can be performed in
different
ways. For example, if $\tau \to 0$ we can get in the limit the $U(\Ga)$
or the skew-commutative algebra with $n$ odd generators of degree 1 and 2n
even generators of degree 2. If $\tau\ne 0$
is a point of second order and $n$
is odd, then $T_n(\CE,\tau)$ degenerates into  polynomials
in $n$ generators of degree 1 and n generators of degree 2.

\noindent
{\it Remark 2.} If $n=3$ algebra $T_3(\CE ,\tau)$ is generated by
three elements of degree one $\{ x_\alpha ,\alpha \in \BZ_3\}$ and
relations can be written as:
\leq
&\theta_0(4\tau)x^2_\alpha+\theta_1(4\tau)x_{\alpha+2}x_{\alpha+1}
+\theta_2(4\tau)x_{\alpha+1}x_{\alpha+2}=C_\alpha, \alpha=0,1,2\cr
&e^{6 \pi i \tau}\theta_\alpha(3\tau)x_\beta C_{\alpha+\beta}
=\theta_\alpha(\tau)\theta_1(2\tau)\theta_2(2\tau)C_{\alpha+\beta}x_\beta
+\theta_{\alpha+1}(\tau)\theta_0(2\tau)\theta_1(2\tau)\cr
&\times C_{\alpha+\beta+1}x_{\beta+1}+\theta_{\alpha+2}(\tau)\cdot
\theta_2(2\tau)\theta_0(2\tau)C_{\alpha+\beta+2}x_{\beta+2}; \alpha,
\beta\in \BZ_3.\cr
\endeq

It is clear that the quotient $T_3(\CE,\tau)/J
\cong Q_3(\CE, 4\tau)$, where $J$ is an ideal generated by $C_\alpha$.
Note also that for generic $\tau$ algebra $T_3(\CE,\tau)$ has two
central elements of degree 3 and 6 and center is generated by these
elements.

\bigskip
\noindent
V. {\it Further Generalizations.}
\par
Let $\Gamma$ be a root system, $\{ \delta_1,\cdots,\delta_h\} \subset
\Gamma^+$-simple positive roots, $(A_{ij})$, $1\le i$,
$j\le h$-the Cartan matrix and
$b_{ij}$-the matrix of scalar products of the simple roots,
$A_{ij}=2b_{ij}/b_{ii}$. We denote by $L$ the lattice with the basis
$\{ \delta_1,\cdots ,\delta_h\}, L\cong \BZ^h$. We will construct
an associative $L$-graded  algebra which is an elliptic deformation
of the universal enveloping algebra
of the following Lie algebra. Let $\GG$
be the
Kac-Moody algebra with the root system $\Gamma, \GG=N_+\oplus f\oplus
N_-$-the Cartan decomposition. Let $\{ g_\gamma, \gamma\in \Gamma^+\}$
be a basis
in $N_+$ and $\varphi_{\gamma_1,\gamma_2}$-the structure constant
$[g_{\gamma_1},g_{\gamma_2}]=\varphi_{\gamma_1,\gamma_2}\cdot g_{\gamma_1
+\gamma_2}$. We denote by $n$ the dominant weight of the $\GG$,
$n_i=\langle n,\delta_i\rangle \ge 0$. Our Lie algebra  $L(\CE ,n)
=L_0\oplus(\oplus L_\gamma), \gamma \in \Gamma_+$. We assign to each
positive root $\gamma$ a line bundle $\xi(\gamma)$ on $\CE$
such that $\xi(\gamma_1)\otimes \xi(\gamma_2)=\xi(\gamma_1+\gamma_2)$ if
$\gamma_1+\gamma_2$ is a root and $\dim H^0(\xi(\delta_i))=n_i$. The
space $L_0$ has a basis $\{ t_1,t_2,\cdots ,t_h\}$ and it is an abelian Lie
algebra; $L_\gamma \cong H^0(\xi(\gamma))$. The bracket
$L_{\gamma_1}\otimes L_{\gamma_2}\to L_{\gamma_1+\gamma_2}$
is: $[f_1,f_2]=g_{\gamma_1,\gamma_2}\cdot f_1\cdot f_2$, where $f_1\cdot f_2$
is just a product of sections. The bracket $[t_i,f_\gamma]={2\over{n_i}}
\langle \delta_i,\gamma\rangle f_\gamma$, where $f_\gamma \in L_\gamma$,
$1\le i\le h$ and $\langle , \rangle$-the scalar product on the root space.
Let $K$ be
the field of meromorphic  functions in variables $t_1,t_2,\cdots ,t_h$
and $U=K\otimes_{\BC[t_1,\cdots ,t_h]}U(L(\CE,n))$. This means that we
add to $U(L(\CE,n))$ the ``arbitrary'' functions in $t_1,\cdots ,t_h$.

Now we construct the $L$-graded algebra $Q_{n,\Gamma}(\CE,\tau)$
which is a graded
deformation of $U$. $Q_{n,\Gamma}(\CE,\tau)=K\oplus(\oplus_lF_l), l$
is the sum
$\sum l_i\delta_i, l_i\ge 0$. The relation between $K$ and $F_l$ is:
$f(t_1,\cdots ,t_h)\cdot u=u\cdot f(t_1+{2\over{n_1}}
\langle \delta_1, l\rangle
\tau, t_2+{2\over{n_2}}\langle \delta_2, l\rangle \tau,\cdots ,t_h
+{2\over{n_h}}\langle \delta_h,l\rangle \tau), u\in F_l$. The space
$F_l$ is a free $K$-module, $F_l=K\tilde \otimes R_l$ and $R_l$ is a
space of meromorphic sections of line bundle $\xi 's^l(\xi)=S^{l_1}\xi
(\delta_1)\boxtimes S^{l_2}\xi(\delta_2)\boxtimes\cdots \boxtimes S^{l_h}\xi
(\delta_h)$ which is a bundle over $S^{l_1}\CE \times S^{l_2}
\CE \times \cdots \times S^{l_h}\CE, \boxtimes$ is outer
tensor product. Elements from $R_l$ we will write as  functions in
$l_1+l_2+\cdots +l_h$ variables: $f(u_{1,1},u_{2,1},\cdots, u_{l_1,1};
u_{1,2},u_{2,2},\cdots ,u_{l_2,2};\cdots ;u_{1,h},\cdots ,u_{l_h,h})$
which is symmetric with respect to the first group of variables,
the second, $\cdots$. Elements from $F_l=K\tilde \otimes R_l$ are functions
in variables $t_1,\cdots ,t_h,u_{11},\cdots ,u_{l_h,h}$ and they are
$\theta$-functions in the variables $\{ u_{ij}\}$ (this is the meaning of
the symbol $\tilde \otimes$). These functions have to satisfy additional
properties.

{\it Description of $F_l$.} (a) function $f$ has a pole of order $\le 1$ on
the divisors $u_{\alpha,j}-u_{\beta,i}-t_j+t_i=0$ if $\langle \delta_i,
\delta_j\rangle \not= 0$ and $i\not= j$, and holomorphic outside the
union of these divisors. (b) functions from $F_l$ are zero on the
following submanifolds: $\{ M_{i,j}\}$. They are labeled by the pairs of
simple roots $\delta_i,\delta_j$, such that $2\langle \delta_i,\delta_j
\rangle/\langle \delta_i,\delta_i\rangle=p_{i,j}<0, p_{ij}\in \BZ$.
Recall, that $p_{ij}$ is an element of Cartan matrix.
The manifold $M_{i,j}$
consist of such $\{ t_\tau,u_{i,j}\}$, that $u_{\alpha_1,i}=u_{\alpha_2,i}
+\langle \delta_i,\delta_i\rangle \cdot \tau, u_{\alpha_2,i}=u_{\alpha_3,i}
+\langle \delta_i,\delta_i\rangle \tau, \cdots u_{\alpha_{-p_{ij}},i}
=u_{\alpha_{-p_{ij}+1}}+\langle \delta_i,\delta_i\rangle \tau,
u_{\alpha_{-p_{ij}+1},i}-t_i=u_{\beta,j}+\langle \delta_i,\delta_j\rangle
\tau-t_j$. A formula for the product in the algebra $Q_{n,\Gamma}
(\CE,\tau)$ is:
$f\in F_l$ and depends on $\{ t_1,\cdots,t_h,u_{i,j}\}, g\in F_l'$ and
depends on $\{ t_1,\cdots,t_h,u'_{i,j}\}, l=\{ l_1,\cdots ,l_h\}, l'=\{ l'_1,
\cdots, l'_h\}$
\leq
&f*g(t_1,\cdots ,t_h,u_{1,1},\cdots,u_{l_1,1},u'_{1,1}\cdots,u'_{l'_1,1};
\cdots ,u_{1,h},\cdots ,u_{l_h,h},u'_{1,h},\cdots ,u'_{l'_h,h})=\cr
&={1\over l_1! l'_1!\cdots l_h! l'_h!}
Symm[\Bigg(
\prod_{{{1\le i,j\le h}\atop{1\le \alpha \le l_i}}
\atop{1\le\beta\le l'_j}}{{\theta(u_{\alpha,i}-u'_{\beta,j}-t_i+t_j-
\langle \delta_i,\delta_j\rangle \tau)}\over
{\theta(u_{\alpha,i}-u'_{\beta,j}-t_i+t_j)}}\Bigg)\times\cr
&
\times f(t_1,\cdots ,t_h,u_{\gamma,\delta})\cdot g(t_1,\cdots ,t_h,
u'_{1,1}-{2\over{n_1}}\langle \delta_i, \delta_1\rangle\tau,\cdots,
u'_{l'_1,1}-{2\over {n_1}}\langle \delta_i,\delta_1\rangle \tau;\cdots\cr
&
u'_{1,h}-{2\over{n_h}}\langle \delta_i,\delta_h\rangle \tau,\cdots,
u'_{l'_h,h}-{2\over{n_h}}
\langle \delta_i,\delta_h\rangle \tau)]\cr
\endeq
The symmetrization here has a usual sense: the expression in
the bracket [ ] has
not the desirable symmetry properties, so we symmetrize to get them.

{\it Main properties} of {\it algebras} $Q_{n,\Gamma}(\CE,\tau)$.

\noindent
(a) If $\Gamma$ is a root system of $A_1$ then $Q_{n,\Gamma}(\CE,
\tau)\cong Q_{n_1}(\CE,{2\over{n_1}}\tau)\bar \otimes K$, where
$K$ is the ring of functions in $t_1$.

\noindent
(b) Suppose that $n=(n_1,\cdots ,n_h)$ is not beside the origin, it means that
all $n_i>1$. Then the algebra $Q_{n,\Gamma}(\CE,\tau)$ is generated by
$ t_1,\cdots ,t_h$ and $F_{\delta_1},F_{\delta_2},\cdots ,
F_{\delta_h}$.

\noindent
(c) Consider an algebra $A(I,\Gamma),I$ is a set of indices; $A(I,\Gamma)$
is generated by $\{ e_{\alpha,i};u_{\alpha,i}\}, \alpha\in I, 1\le i\le h$,
the relations are:
\leq
e_{\alpha,i}e_{\beta,j}&=-e^{2\pi i(u_{\beta,j}-u_{\alpha,i}-t_j+t_i)}
{{\theta(u_{\alpha,i}-u_{\beta,j}-t_i+t_j-\langle \delta_i,
\delta_j\rangle \tau)}
\over{\theta(u_{\beta,j}-u_{\alpha,i}-t_j+t_i-\langle \delta_i,\delta_j
\rangle \tau)}}e_{\beta,j}\cdot e_{\alpha,i}\cr
e_{\alpha,i}u_{\beta,j}&=(u_{\beta,j}-{2\over{n_j}}\langle \delta_i,
\delta_j\rangle \tau)e_{\al,i};\hbox{here}\ \al\not= \beta\ \hbox{when}\ i=j\cr
e_{\al,i}u_{\al,i}&=(u_{\al,i}+(2-{4\over{n_i}})\tau)e_{\al,i};
e_{\al,i}t_j=(t_j-{2\over{n_j}}\langle \delta_i,\delta_j\rangle \tau)e_{\al,i}
\cr
\endeq
There is a homomorphism $Q_{n,\Gamma}(\CE,\tau)\to A(I,\Gamma), f\in
F_{\delta_i}$
goes to $\sum_\al f(u_{\al,i})e_{\al,i}$

\bigskip
\noindent
VI. Serre relations and classical limits.

Algebra $Q_{n,\Gamma}(\CE,\tau)$ is generated by $\{ t_i\}$ and
$\{ F_{\delta_i}\}$ only in the case when $\tau$ is not a point of
finite order. If $\tau$ is a point of finite order, the subalgebra
$\tilde Q_{n,\Gamma}(\CE,\tau)$ generated by $\{ t_i\}$ and
$\{ f_{\delta_i}\}$ is smaller. For example, suppose $\tau^N=0$ on $\CE$,
$\tau\not= 0$. In this case functions from $\tilde Q_{u,\Gamma}(\CE,\tau)$
have additional zeros. Namely, $f(t_1,\cdots,t_h;u_{\gamma,\delta})$ is zero
if $u_{\al_1,i}=u_{\al_2,i}+\langle \delta_i,\delta_i\rangle \tau,
u_{\al_2,i}=u_{\al_3,i}+\langle \delta_i,\delta_i\rangle \tau,\cdots ,
u_{\al_{N-1},i}=u_{\al_N,i}+\langle \delta_i,\delta_i\rangle \tau$ for
each $i, \al_1,\cdots ,\al_N$. It is natural to conjecture that
$\tilde Q_{n,\Gamma}(\CE,\tau)$ coincides with the space of functions
with this property. All these zero conditions (see also point (b) from
the description of $F_l$) are a sort of Serre relations in the universal
enveloping of the nilpotent part in universal enveloping of
$q$-deformed Kac-Moody.

To make it clear consider the following construction. Let $L$ be a lattice
and $D=\oplus D_l, l\in L$ is $L$-graded associative algebra;
$\delta_1,\cdots ,\delta_h$-basis in $L$ and the space $D_l, l=l_1\delta_1
+l_2\delta_2+\cdots +l_k\delta_k$ is a space of functions $f$ in
variables $u_{1,1},u_{2,1},\cdots ,u_{l_1,1};u_{1,2},u_{2,2},
\cdots ,\hfil\break\noindent
u_{l_2,2};\cdots;u_{1,h},u_{2,h},\cdots ,u_{l_h,h};$
which are symmetric with respect to the first group of $l_1$ variables,
with resepct to the second group of $l_2$ variables, $\cdots$. The product
is given by our standard formula:
\leq
&f*g(u_{1,1},\cdots ,u_{l_1,1},u'_{1,1},\cdots ,u'_{l'_1,1};\cdots ;
u_{1,h},\cdots,u_{l_h,h},u'_{1,h},\cdots u'_{l'_h,h})\cr
&=Symm[\bigg(
\prod_{{{1\le i,j\le h}\atop{1\le \al\le l_i}}\atop{1\le\beta
\le l'_j}}\lambda_{i,j}(u_{\al,i},u'_{\beta,j})\bigg)
\times f(u_{\gamma,\delta})
\cdot g(u'_{\gamma,\delta})]\cr
\endeq
Here $f\in D_l, l=\sum l_i\delta_i$ and depends on $U_{\gamma,\delta}$ and
$g\in D_{l'}, l'=\sum l'_i\cdot \delta_i$ and depends on $u'_{\gamma,\delta}$
and $f\times g\in D_{l+l'}$ and depends on $u_{\gamma,\delta},
u'_{\gamma,\delta}; \{ \lambda_{i,j}\}$ is an arbitrary set of function
in two variables. Suppose for simplicity that the functions $\lambda_{i,j}$
are regular and let $K_{i,j}$ be a divisor of zeros of the function
$\lambda_{i,j}$. Let us denote by $\tilde D$ the subalgebra in $D$ generated
by $\oplus_i D_{\delta_i}$.

\proclaim
Proposition. Functions from $\tilde D$ have zeros on the set of
$\{ u_{\gamma,\delta}\}$ such that there
exists the subsequence of points:
$\{ u_{1,i_1},u_{2,i_1},\cdots ,u_{j_1,i_1}; u_{1,i_2},u_{2,i_2},\cdots ,
u_{j_2,i_2};\cdots ;u_{1,i_s},\cdots u_{i_s,i_s}\}$ such that
$\{ u_{1,i_1},u_{2,i_1}\}\in K_{i_1,i_1};\{ u_{2,i_1},u_{3,i_1}\}
\in K_{i_1,i_1};\cdots ;\langle u_{j_1-1,i_1,}u_{j_s,i_s}\}\in
K_{i_1,i_1},\{ u_{i_1,i_1},u_{1,i_2}\}\in K_{i_1,i_2},\{ u_{1,i_2,}u_{2,i_2}\}
\in K_{i_2,i_2},\cdots ,\{ u_{j_s-1,i_s},u_{i_1,i_s}\}\in K_{i_s,i_s},$
and the last is $\{ u_{j_s,i_s},u_{1,i_1}\}\in K_{i_s,i_1}$. Algebra
$Q_{n,\Gamma}(\CE,\tau)$ is a subalgebra in the twisted version of $D$
with some special $\{ \lambda _{i,j}\}$ they have poles but it is not
essential, the situation is the same. It is interesting that in this case
these zero conditions are consequences of the conditions for
$\{ u_{1,i_1},u_{2,i_1},\cdots,u_{j_1,i_1};u_{1,i_2}\}$.

When $\tau\to 0$
algebra $Q_{n,\Gamma}(\CE,\tau)$ tends to commutative algebra. Therefore
it is possible to define the quasi-classical limit of
$Q_{n,\Gamma}(\CE,\tau)$, to find out symplectic leaves and so on.
Another way to construct the $\tau\to 0$ limit is the following.
Using the considerations of the point I we first have the braided tensor
product $Q_{n,\Gamma}\bar \otimes Q_{n,\Gamma}$ and the comultiplication
map $Q_{n,\Gamma}\to Q_{n,\Gamma}\bar \otimes Q_{n,\Gamma}$. The simplest
definition of $Q_{n,\Gamma}\bar \otimes Q_{n,\Gamma}$ can be done using
algebra $A(I,\Gamma)$ from point V. Suppose I is infinite, then
$Q_{n,\Gamma}(\CE,\tau)\to A(I,\Gamma)$ is imbedding. Decompose I into
two infinite subsets $I_1\cap I_2=I$, so $A(I_j,\Gamma)\subset A(I,\Gamma)$
and two copies of $Q_{n,\Gamma}$ is contained in $A(I,\Gamma)$: one
is $Q^{(1)}_{n,\Gamma}\subset A(I_1,\Gamma)$ and the second
$Q^{(2)}_{n,\Gamma}\subset A(I_2,\Gamma)$. These $Q^{(1)}_{n,\Gamma}$
and $Q^{(2)}_{n,\Gamma}$ generate $Q_{n,\Gamma}\bar \otimes Q_{n,\Gamma}.$
Comultiplication is given by the familiar formula
$f\in F_{\delta_i}\to f\otimes 1+1\otimes f, t_j\to t_j\otimes 1+1
\otimes t_j$. When $\tau \to 0$ comultiplication degenerate into the
non-trivial operator $Q_{n,\Gamma}(\CE,0)\to Q_{n,\Gamma}(\CE,0)\otimes
Q_{n,\Gamma}(\CE,0)$. Hence we get a Hopf algebra. Dual object is algebra
$U(L,\CE)$-the universal enveloping of the Lie algebra from the
beginning of point V. Informally it means that algebra $Q_{n,\Gamma}(\CE,0)$
is more or less an algebra of functions on a group of currents from
$\CE$ into maximal nilpotent subgroup into Kac-Moody algebra.
Algebra $Q_{n,\Gamma}(\CE,\tau)$ for generic $\tau$ is generated
by $\{ t_i\}$ and the sum of $\{ F_{\delta_i}\}$.
When $\tau \to 0$, then the relations between $\{ F_{\delta_i}\}$ tend
to Serre relations in $U(L,\CE)$.

\remark There is a way to generalize our construction. Algebra from the point
IV is connected with maximal nilpotent subalgebra in Lie superalgebra
osp(1,2) which has a basis $l_1,l_2$ and a bracket $[l_1,l_2]=l_2,
[l_2,l_1]=0, \deg l_1=1, \deg l_2=2.$ Changing the functions $\lambda_{i,j}$
we can cover many superalgebras. It will be a subject of another paper.
\par

\beginsection
\S3

Let $N_{n,k}$
be the following moduli space.
Consider the moduli space of indecomposible $k$-dimensional
bundles of degree $n$.
More presicely we need corresponding universal family
of bundles
$\xi_{n,k}(z)$,
$z\in\BC\to \BC/\Gamma$,
is a parameter on the elliptic curve
$\BC/\Gamma=\CE$.
Element from the space
$N_{n,k}$
is a pair
$(z,\beta)$, $z\in\CE$
and $\beta\in \Ext^1(\xi_{n,k}(z),\xi_{0,1})$.
Define an algebra
$\tilde Q_{n,k}(\CE,\tau)$
with generators
$x_1,\ldots, x_n$, $z$.
Elements $\{x_i\}$
obey the same relations as the generators of algebra
$ Q_{n,k}(\CE,\tau)$
(see the begining of introduction),
$x_iz=(z+\tau')x_i$,
$\tau^1\in\BC$.
Algebra
$\tilde Q_{n,k}(\CE,\tau)$
consists of expressions
$\sum x^{i_1}_1\ldots x^{i_n}_n\cdot
f_{i_1,\ldots, i_n}(z)$,
where $f$ are elliptic functions associated with
$\CE$.
If $\tau$ and $\tau '$ go to zero
$(\tau/\tau'$ is generic) we get the hamiltonian structure
on $N_{n,k}$.
\par
Let \Mod\ be a set of vector bundles on
$\CE$ defined up to isomorphism.
Define an action of $\CE$ on \Mod\ by the formula:
$\beta\in\CE$ acts as
\eq
&T_\beta : \xi_{n_1,k_1}(\alpha_1)
\oplus \xi_{n_2,k_2}(\alpha_2)
\oplus \cdots
\oplus \xi_{n_p,k_p}(\alpha_p)
\to
\xi_{n_1,k_1}(\overline\alpha_1)
\oplus \cdots
\oplus \xi_{n_p,k_p}(\overline\alpha_p)
\cr
\endeq
where
$\overline\alpha_i=\alpha_i+\beta\cdot
(n_i(k+1)-k_i\cdot n)$,
$n=n_1+n_2+\cdots +n_p$,
$k+1=k_1+\cdots +k_p$.
\par
There is a map
\eq
&\theta : N_{n,k}\to \Mod \to \Mod/\CE.
\cr
\endeq
Here $\Mod/\CE$
is a space of orbits of $\CE$ acting on
\Mod\ by $\beta\to T_\beta$.
\proclaim
Proposition.
Symplectic leaves of the hamiltonian structure on
$N_{n,k}$ are the fibers of map $\theta:
N_{n,k}\to \Mod/\CE$.
\par
Now consider the case of algebra
$Q_{n,k}(\CE,\tau)$
without the additional variable $z$.
In the introduction we described the symplectic leaves
of the hamiltonian structure on $\BC P^{n-1}$.
This $\BC P^{n-1}$
is a moduli space of exact sequences $0\to\xi_{0,1}\to ? \to
\xi_{n,k}(z)\to 0$
where $z$ is fixed.
There is a birational map
\eq
&\BC P^{n-1}\to \Mod^s_{n,k}(z).\cr
\endeq
Here $\Mod^s_{n,k}(z)$
is the moduli space of stable $k+1$-dimensional
bundles of degree $n$ with fixed determinant.
The determinant have to be equal to the determinant of
$\xi_{0,1}\oplus \xi_{n,k}(z)$.
Due to \S1 we know that
$\Mod^s_{n,k}(z)$
is isomorphic to
$\BC P^{c-1}$,
where $c$ is the maximal common divisor of
$n$ and $k+1$.
It is possible to draw the following commutative diagram:
$$
\matrix{
\BC^n & \smash{\mathop{\longrightarrow}\limits^{\bar\theta}}
& \BC^c \cr
\llap{$\vcenter{\hbox{$\scriptstyle P_1$}}$}\Big \downarrow
& & \Big \downarrow \rlap {$\vcenter{\hbox{$\scriptstyle P_2$}}$} \cr
\BC P^{n-1} & \longrightarrow
& \Mod_{n,k}(z)\cong \BC P^{c-1} \cr}
$$
Here $P_1$ and $P_2$
are natural projections.
They are not defined in the origins.
The map $\bar\theta$ is regular and can be written as
a set $\{\theta_1,\ldots,\theta_c\}$
of polynomials of degree $n/c$.
It is reasonable to think about $\BC^c\backslash \{0\}$
as a space of determinant bundle over
$\Mod_{n,k}(z)$.
So, the meaning of our diagram is that the rational map
in the bottom is covered by the regular map of determinant
bundles.
\proclaim
Propositions.
The symplectic leaves of hamiltonian structure on
$\BC^n$ (it corresponds to the classical limit of
$Q_{n,k}(\CE,\tau)$)
are the intersections of the fibers of the map
$\BC^n\to \Mod \to \Mod/\CE$
and $\BC^n\to \BC^c$.
In particular, $\bar\theta$
is hamiltonian map, if the hamiltonian
structure on $\BC^c$ is trivial.
\par
Using this proposition and results from \S1 we can describe all
symplectic leaves of the classical limit of
$Q_{n,k}(\CE,\tau)$.
\par
The similar results are true in the case of Borel bundles
on $\CE$. Consider for example the $sl_r$-case.
So, let $M_{n_1,\ldots,n_{r-1}}$
be a moduli space of objects:
$\Gamma$-dimensional bundle $\nu$ with filtration
$\nu_0\subset \nu_1\subset \cdots \subset \nu_{r-1}$
and the isomorphisms:
$\nu_0\cong \xi_{0,1}(z_1)$,
$\nu_1/\nu_0\cong \xi_{n_1,1}(z_2),\ldots,
\nu_{r-1}/\nu_{r-2}\cong \xi_{n_{r-1},1}(z_r)$.
Suppose also that the determinant of $\nu$ is fixed.
The classical limit of the algebras
$Q_{n,\Gamma}(\CE,\tau)$,
where $n=(n_1,n_2,\ldots, n_{r-1})$
and $\Gamma$ is a root system of $sl_r$
gives us a hamiltonian structure on
$M_{n_1,\ldots,n_{r-1}}$.
As in
$Q_{n,k}(\CE,\tau)$
case it is a hamiltonian map
\eq
&M_{n_1,\ldots,n_{r-1}}\to \BC^h.\cr
\endeq
Here $\BC^h$ has a trivial hamiltonian structure and $h$ is a maximal common
divisor of $r$ and $n_1+2n_2+\cdots +(r-1)n_{r-1}$.
It is also possible to define a map
$M_{n_1,\ldots,n_{r-1}}\to \Mod\to \Mod/\CE$
and description of symplectic leaves is the same as in
$Q_{n,k}$-case.

\beginsection
References

\item{1.}
M.F. Atiyah, Vector bundles over an elliptic curve,
Proc.\ London Math.\ Soc.\ 7 (1957), 414-452.

\item{2.}
I. Grojnowski, Delocalized equivariant elliptic cohomology,
Yale University preprint, 1994.

\item{3.}
A.V. Odesskii and B.L. Feigin, Sklyanin Algebras
Associated with an Elliptic Curve [in Russian],
Institute for Theoretical Physics, Kiev (1988).

\item{4.}
A.V. Odesskii and B.L. Feigin, ``Sklyanin's elliptic
algebras'', Funkts.\ Anal.\ Philozhen., 23, No.3, 45-54 (1989).

\item{5.}
A.V. Odesskii and B.L. Feigin, ``Constructions of Sklyanin elliptic
algebras and quantum $R$-matrices'',
Funkts.\ Anal.\ Philozhen., 27, No.1, 37-45 (1993).

\item{6.}
A.V. Odesskii, Rational Degeneration of Elliptic Quadratic
Algebras, RIMS 91, Project ``Infinite Analysis'',
June 1 - August 31, 1991.

\item{7.}
A.V. Odesskii and B.L. Feigin, Sklyanin Elliptic
Algebras. The case of points of finite order,
RIMS-986, July 1994, Kyoto University, Kyoto, Japan.

\end